\documentclass[3p]{elsarticle}
\usepackage{bm}
\usepackage[parfill]{parskip}
\usepackage[utf8]{inputenc}
\usepackage{url}
\usepackage{amsmath,amssymb,amsfonts,amsthm}
\usepackage{color}
\usepackage{bm}
\usepackage{graphicx}
\usepackage{subcaption}
\usepackage{lipsum}
\usepackage{amsfonts}
\usepackage{graphicx}
\usepackage{epstopdf}
\usepackage{array}
\usepackage[hidelinks]{hyperref}
\usepackage{lineno}
\usepackage{bbm}
\usepackage{xfrac}
\usepackage{xcolor}

\usepackage[capitalize,noabbrev]{cleveref}
\usepackage{amssymb}
\usepackage{pifont}
\usepackage{amsmath}
\usepackage{algorithm}
\usepackage[noend]{algpseudocode}
\usepackage{booktabs}       
\usepackage[textsize=tiny]{todonotes}
\newcommand{\cmark}{\ding{51}}%
\newcommand{\xmark}{\ding{55}}%

\newcommand{\green}[1]{{#1}}

\ifpdf
  \DeclareGraphicsExtensions{.eps,.pdf,.png,.jpg}
\else
  \DeclareGraphicsExtensions{.eps}
\fi

\makeatletter
\newsavebox\myboxA
\newsavebox\myboxB
\newlength\mylenA

\newcommand*\xoverline[2][0.75]{%
    \sbox{\myboxA}{$\m@th#2$}%
    \setbox\myboxB\null
    \ht\myboxB=\ht\myboxA%
    \dp\myboxB=\dp\myboxA%
    \wd\myboxB=#1\wd\myboxA
    \sbox\myboxB{$\m@th\overline{\copy\myboxB}$}
    \setlength\mylenA{\the\wd\myboxA}
    \addtolength\mylenA{-\the\wd\myboxB}%
    \ifdim\wd\myboxB<\wd\myboxA%
       \rlap{\hskip 0.5\mylenA\usebox\myboxB}{\usebox\myboxA}%
    \else
        \hskip -0.5\mylenA\rlap{\usebox\myboxA}{\hskip 0.5\mylenA\usebox\myboxB}%
    \fi}
\makeatother

\journal{arXiv.org}

\newcommand{\TheTitle}{Construction of high-order conservative basis-update \& Galerkin dynamical low-rank integrators\footnote{This manuscript has been authored by UT-Battelle, LLC under Contract No. DE-AC05-00OR22725 with the U.S. Department of Energy. The United States Government retains and the publisher, by accepting the article for publication, acknowledges that the United States Government retains a non-exclusive, paid-up, irrevocable, world-wide license to publish or reproduce the published form of this manuscript, or allow others to do so, for United States Government purposes. The Department of Energy will provide public access to these results of federally sponsored research in accordance with the DOE Public Access Plan(\url{http://energy.gov/downloads/doe-public-access-plan}).}\footnote{The experimental results can be reproduced with using \url{https://github.com/ScSteffen/Publication_High_Order_Conservative_Integrators}.}} 

\date{\today}

\usepackage{amsopn}

\journal{arXiv}

\author[adressInnsbruck]{Lukas Einkemmer}
\author[adressAs]{Jonas Kusch}
\author[adressKarlsruhe]{Steffen Schotth\"ofer}

\address[adressInnsbruck]{University of Innsbruck, Numerical Analysis and Scientific Computing, Innsbruck, Austria}
\address[adressAs]{Norwegian University of Life Sciences, Scientific Computing, \r{A}s, Norway}
\address[adressKarlsruhe]{  Computer Science and Mathematics Division,
  Oak Ridge National Laboratory,
  Oak Ridge, TN 37831 USA}

\begin{document}
\begin{frontmatter}

\title{\TheTitle}


\begin{abstract}
Numerical simulations of kinetic problems can become prohibitively expensive due to their large memory requirements and computational costs. A method that has proven to successfully reduce these costs is the dynamical low-rank approximation (DLRA). A major accomplishment in the field of DLRA has been the derivation of robust time integrators that are not limited by the stiffness of the DLRA evolution equations. One key question is whether such robust time integrators can be made locally conservative, i.e., can they preserve the invariants and associated conservation laws of the original problem? In this work, we propose extensions to commonly used basis-update \& Galerkin (BUG) integrators that preserve invariants of the solution as well as the associated conservation laws with little or no additional computational cost. This approach requires only minor modifications of existing implementations. The properties of these integrators are investigated by performing numerical simulations in radiation transport and plasma physics.
\end{abstract}

\begin{keyword}
Dynamical low-rank approximation, structure-preserving integrator, kinetic equations, basis-update \& Galerkin integrator
\end{keyword}

\end{frontmatter}

\section{Introduction}

A significant number of problems in physics and engineering require the solution of high-dimensional partial differential equations. Their numerical discretization is extremely expensive in terms of both memory and computational cost. A promising model order reduction method that can significantly reduce computational costs and memory requirements, while maintaining good accuracy, is the dynamical low-rank approximation (DLRA), which has been proposed in \cite{koch2007dynamical}. The main idea of DLRA is to restrict the evolution to the manifold of low-rank functions. In this case, the stability region of classical time integration schemes depends on the curvature of the manifold which can significantly reduce allowable time step sizes. The perhaps most frequently used integrators that do not depend on this curvature are the projector--splitting integrator \cite{lubich2014projector,kieri2016discretized} and the basis-update \& Galerkin (BUG) integrator \cite{ceruti2022unconventional,ceruti2022rank}. The latter has been extended to perform time updates of all low-rank factors in parallel, resulting in a scheme known as the parallel BUG integrator \cite{ceruti2023parallel}. Extensions of these first-order integrators to higher order have been proposed in \cite{ceruti2024robust} and \cite{kusch2024second}.
The development of robust time integrators for DLRA has enabled the use of dynamical low-rank approximation for a wide variety of problems. Following its first use for kinetic problems \cite{einkemmer2018low}, DLRA has gained a considerable amount of attention in this research field \cite{einkemmer2024review}.

One particularly important question is whether DLRA is able to conserve the physically relevant invariants and their corresponding conservation laws. For problems such as the Schrödinger equation, this is immediate (as the $L^2$ projection preserves the corresponding norm); in a similar way, gradient flow and loss-descent preservation (up to truncation tolerance) can be shown in DLRA for optimization \cite{ceruti2022rank,NEURIPS2022_7e98b00e}. There are also some results on Hamiltonian systems \cite{ceruti2022rank}. However, for problems with an underlying local conservation law, in particular, kinetic equations, it is well known that applying DLRA results in a set of equations that, even if solved exactly (i.e.~without space and time discretization), do violate the underlying physical structure. The first step to remedy this was taken in \cite{einkemmer2021mass}, where the Galerkin condition was modified to result in a set of equations for the low-rank factors that preserve mass, momentum, and energy for Vlasov--Poisson systems. It has also been shown that this can be combined with time and space discretization strategies that result in a fully discrete conservative method. However, the approach presented in this work is not robust. The main issue is that the equations for the low-rank factors are modified in such a way that they lose the property required for applying the projector splitting integrator. Realizing this, in \cite{einkemmer2023robust}, the recently proposed augmented BUG integrator \cite{ceruti2022rank} has been used to integrate the equations of motion forward in time. This results in a conservative scheme up to the truncation step. However, in the same work, it is shown that the truncation can be performed in a conservative way (a similiar approach is also used in \cite{guo2022conservative}, where a conservative truncation is needed in the context of low-rank step and truncate methods). This yields a scheme that is both conservative and robust, see also \cite{uschmajew2025discontinuous} for an extension to discontinuous--Galerkin discretizations. There are two main disadvantages of this approach. First, the scheme requires the augmentation with additional basis functions, which increases computational cost. Second, the scheme is limited to a forward Euler time integration scheme in the coefficient update. This is problematic since it restricts accuracy and limits the use of stable time integration methods.

An alternative to constructing locally conservative integrators is a micro-macro decomposition, as was proposed in \cite{koellermeier2024macro} for a shallow water moment equation. Such conservative micro-macro decompositions have also been suggested for nodal discretizations for Vlasov--Poisson in \cite{coughlin2024robust} and thermal radiative transfer in \cite{frank2025asymptotic}. Further methods that are similar to a macro-micro decomposition, which do, however, not strictly fall into the category of the works mentioned before include \cite{baumann2024energy} and \cite{patwardhan2025parallel}. Such approaches are effective to preserve invariants and allow for a combination with higher-order time integration methods. They, however, commonly require bigger modifications to a given code base and can require additional computational costs to ensure the micro component of the solution does not impact conservation. Further conservative low-rank approaches that do not fall into the realm of DLRA are conservative variants of step and truncate methods, see, e.g., \cite{allmann2022parallel,guo2024local1,guo2024local2}.

In this paper, we propose a different approach. We investigate local conservation on a time-discretized level, where we can show that local conservation is preserved if an (at most small) additional augmentation is made. This allows us to modify high-order BUG integrators such as \cite{ceruti2024robust} to allow for conservation when the DLRA evolution equations are solved with arbitrary explicit Runge--Kutta methods. Based on this discussion, we can also show that the augmented BUG integrator \cite{ceruti2022rank} is, by construction, conservative when forward Euler time updates are used in the $K$, and $S$-substeps and when combined with a conservative truncation. More precisely, when directly applied to the classic equations of motion for the low-rank factors (i.e.~without ensuring that the underlying equations are conservative, as in \cite{einkemmer2023robust}), the integrator is conservative. Lastly, we propose a minor modification to the parallel integrators \cite{ceruti2023parallel,kusch2024second} which makes them locally conservative without requiring additional augmentation. The main idea is to replace the coefficient update with a $K$-step, which reduces computational costs and directly yields a conservative method.

Let us emphasize that when we consider conservation in this paper, we do not only consider (global) conservation of invariants (this can be done in a relatively straightforward way by introducing Lagrange multipliers; see, e.g., \cite{einkemmer2019quasi,peng2021high}). Instead, our goal is to obtain numerical methods that also preserve the corresponding (local) conservation laws. In the PDE case, the latter is more important in order to improve the quality of the solution \cite{filbet2001conservative,perse2021geometric,taitano2021conservative,einkemmer2019quasi}. We also note that all the methods considered in the paper are based on the augmented BUG integrator and are thus naturally rank adaptive~\cite{ceruti2022rank}.

The paper is structured in the following way: after the introduction we provide a recap on kinetic equations and their conserved quantities in Section~\ref{sec:kin}. Dynamical low-rank approximation is reviewed in Section~\ref{sec:DLRA_main}, where special focus is put on the derivation of the DLRA evolution equations in Section~\ref{sec:DLRA} followed by a discussion of the {augmented} BUG integrator in Section~\ref{sec:rBUG}. Section~\ref{sec:cBUG} provides a short recap on continuously conservative and macro-micro methods. In Section~\ref{sec:local-cons-rbug} we discuss the ingredients for deriving locally conservative BUG integrators. Section~\ref{sec:examples} discusses two special cases in which a conservative truncation is already sufficient to obtain a locally conservative integrator. In Section~\ref{sec:cons-parallel}, a different approach is discussed that holds for parallel BUG integrators and enables a computationally inexpensive modification to achieve local conservation. The conservation and approximation properties of the locally conservative BUG integrators are shown and compared for radiation transport in Section~\ref{sec:results_rad} and plasma physics in Section~\ref{sec:vp-numerics}. In the case of the Vlasov--Poisson equations, the approximation of the electric field, the choice of a weight function in velocity, and mitigating effects of rounding errors are crucial --- not only in this work, but for conservation in general. We dedicate Section~\ref{sec:discussion_vlasov} to provide a detailed discussion of these points and how to robustly and reliably achieve conservation in this setting.

\section{Recap: Kinetic equations}\label{sec:kin}

Kinetic equations govern the physical processes in various applications. Commonly, they are based on a mesoscopic description of a large number of particles defined in a spatial domain $D_x\subset \mathbb{R}^3$ moving with a velocity $v\in D_v\subset \mathbb{R}^3$. The state of the system is described by a particle density $f : \mathbb{R}_+ \times D_x \times D_v \rightarrow \mathbb{R}$, which is evolved according to a kinetic evolution equation
\begin{align*}
    \partial_t f(t,x,v) = F(f(t,x,v)),
\end{align*}
where $F$ usually models streaming (particles moving in direction $v$) and interaction (particles collide with a background material or react to an electric field). A key property of such systems is that they obey local conservation laws. That is, given a vector of conserved basis functions $U(v)$, the dynamics fulfills a local conservation law for the quantities $\Phi(t,x) = \int_{D_v} f(t,x,v)\,dv$, i.e.,
\begin{align*}
    \partial_t \Phi(t,x) = \int_{D_v} F(f(t,x,v)) U(v)\,dv, \quad \text{ where } \int_{D_x}\int_{D_v} F(f(t,x,v)) U(v)\,dv dx = 0.
\end{align*}
The latter holds if no inflow/outflow exist at the spatial boundary and, in this case, ensures global conservation, i.e., $\frac{d}{dt}\int_{D_x}\Phi(t,x)\,dx = 0$. In the following, we discuss kinetic equations in two specific application fields, namely radiation transport and plasma physics, in more detail.

\subsection{Radiation transport}

In the case of radiation transport, radiation particles move with a fixed velocity, which can be scaled to unit velocities, i.e., $v\in\mathbb{S}^2$. The dynamics of the particle density $f$ is then described by the radiative transfer equation
\begin{align}\label{eq:RTE}
    \partial_t f(t,x,v) + v \cdot \nabla_x f(t,x,v) + \sigma_a(x) f(t,x,v) = \int_{\mathbb{S}^2} k(x, v, v') (f(t,x,v')-f(t,x,v))\, dv',
\end{align}
where $\sigma_a$ is the absorption cross-section and $k(x, v, v') = k(x, v', v)$ is the symmetric scattering kernel. 
Defining the scalar flux $\phi(t,x) := \int_{\mathbb{S}^2}f(t,x,v)\,dv$, integrating \eqref{eq:RTE} with respect to $v$ gives the local conservation law
\begin{align*}
    \partial_t\phi(t,x) + \sigma_a(x)\phi(t,x)  + \nabla_x \cdot \int_{\mathbb{S}^2}v f(t,x,v)\, dv = 0.
\end{align*}
If no absorption is present and the solution is zero at the boundary $\partial D_x$ or periodic boundary conditions are used, integrating the local conservation law over time yields $\frac{d}{dt}\int_{D_x}\phi(t,x)\, dx = 0$, i.e., the total mass is conserved.

\subsection{Vlasov equation}
Another prominent example of a kinetic equation is the Vlasov equation. Assuming a constant background ion density, the Vlasov equation determines in electron density $f(t,x,v)$ for a collisionless plasma in the electrostatic regime according to
\begin{align*}
    \partial_t f(t,x,v) + v\cdot \nabla_x f(t,x,v) - E(f)(t,x) \cdot \nabla_v f(t,x,v) =&\, 0\;,\\
    \nabla_x E(f)(t,x) = 1 - \int_{\mathbb{R}^3} f(t,x,v)\,dv, \quad \nabla_x \times E(f)(t,x) =&\, 0\;,
\end{align*}
where for dimension $d$, we have $x\in D_x \subset \mathbb{R}^d$ and $v\in D_v \subset \mathbb{R}^d$. Commonly, one is interested in conserving mass $M$, momentum $J$ and energy $\mathcal{E}$, given by
\begin{align*}
    &M(t) = \int_{D} f(t,x,v)\, dxdv, \quad J(t) = \int_{D} v f(t,x,v)\,dx dv,\\
    \text{ and } &\mathcal{E}(t) = \frac12 \int_D v^2 f(t,x,v)\,dxdv + \frac12 \int_{D_x} E(t,x)^2\,dx\;,
\end{align*}
where $D = D_x \times D_v$ and $v^2$, $E^2$ denote squared Euclidean norms, e.g. $v^2 = \sum_{i=1}^d v_i^2$. The Vlasov equation then fulfills a local conservation law for the local mass $\rho(x)$, momentum $j(x)$, and energy $e(x)$ which are defined as
\begin{align*}
    &\rho(t,x) = \int_{D_v} f(t,x,v)\, dv, \quad j(t,x) = \int_{D_v} v f(t,x,v)\, dv,\\
    \text{ and } &e(t,x) = \frac12 \int_{D_v} v^2 f(t,x,v)\,dv + \frac12 E(t,x)^2\;.
\end{align*}
The local conservation laws for these quantities then read
\begin{align*}
    \partial_t \rho(t,x) + \nabla_x \cdot j(t,x) =&\, 0 \\
    \partial_t j(t,x) + \nabla \cdot \int_{D_v} (v \otimes v) f(t,x,v)\, dv =&\, - E(t,x)\rho(t,x) \\
    \partial_t e(t,x) + \nabla_x \cdot \frac12 \int_{D_v} vv^2 f(t,x,v)\, dv =&\, E(t,x) \cdot (\partial_t E(t,x)- j(t,x))\;.
\end{align*}
Conservation of mass, momentum and energy is guaranteed since $E(1-\rho) = \nabla_x\cdot (E\otimes E - \frac12 E^2)$ and $\int_{D_x} E \,dx = 0$.

In the following, to make the notation more general, we denote our kinetic equation as $\partial_t f(t,x,v) = F(f(t,x,v))$. Given a vector of $m$ conservative basis functions $U(v)\in\mathbb{R}^m$ and a vector of conserved quantities $\Phi(t,x):= \int_{D_v} f(t,x,v) U(v)\, dv$, we wish to investigate dynamical low-rank time integrators which preserve the local conservation law
\begin{align*}
    \partial_t\Phi(t,x) + \int_{\mathbb{S}^2}F(f(t,x,v)) U(v)\, dv = 0,
\end{align*}
which when integrating over the spatial domain directly yields conservation of $\int_{D_x}\Phi(t,x)\, dx$ for periodic or Dirichlet boundary conditions $f\big|_{\partial D_x} = 0$.

\section{Recap: Dynamical low-rank approximation}\label{sec:DLRA_main}

\subsection{Low-rank evolution equations}\label{sec:DLRA}
Dynamical low-rank approximation for kinetic problems revolves around the fundamental concept of expressing and evolving the solution on a manifold consisting of functions with a rank of $r$, which we denote as
\begin{align*}
    \mathcal{M}_r = \Bigg\{&f\in L^2(D_x\times D_v) : f(x,v) = \sum_{i,j=1}^r X_i(x)S_{ij} V_{j}(v) \text{ with invertible } S = (S_{ij})\in\mathbb{R}^{r\times r},\\ 
    &X_i \in L^2(D_x), V_j\in L^2(D_v) \text{ and } \langle X_i, X_j\rangle_x = \delta_{ij}, \langle V_i, V_j\rangle_{v} = \delta_{ij} \Bigg\}.
\end{align*}
Here, we denote the $L^2$ inner product\footnote{Note that when $D_x,D_v = \mathbb{R}^d$, the inner products $\langle \,\cdot\,, \,\cdot\,\rangle_{x}$ and $\langle \,\cdot\,, \,\cdot\,\rangle_{v}$ can be defined as a weighted inner product to ensure integrability of conservative basis functions $U(v)$.} over the spatial and velocity domain as $\langle \,\cdot\,, \,\cdot\,\rangle_{x}$ and $\langle \,\cdot\,, \,\cdot\,\rangle_{v}$ respectively. That is, we represent the solution through $r$ orthonormal basis functions $X_i$ and $V_i$ in space and velocity respectively. The corresponding expansion coefficients are collected in the coefficient matrix $S$. Then, at a given time $t$, a low-rank solution can be written as
\begin{align*}
    f(t,x,v) = \sum_{i,j = 1}^r X_i(t,x)S_{ij}(t) V_{j}(t,v) := X(t,x)^{\top}S(t) V(t,v).
\end{align*} Here, for efficiency of notation we define $X(t,x) = (X_i(t,x))_{i=1}^r \in\mathbb{R}^r$ and $V(t,v) = (V_i(t,v))_{i=1}^r \in\mathbb{R}^r$. To confine the dynamics of the solution within the low-rank manifold $\mathcal{M}_r$, it is necessary to identify the associated tangent space at a given position $f$. To ensure a bijective map between the manifold and its tangent space, we impose Gauge conditions $\langle \dot X_i, X_j\rangle_x = 0$ and $\langle \dot V_i, V_j\rangle_{v} = 0$. Then, following \cite{koch2007dynamical,einkemmer2018low}, the tangent space can be written as
\begin{align*}
    \mathcal{T}_f\mathcal{M}_r = \Bigg\{&\delta f\in L^2(D_x\times D_v) : \delta f(x,v) = \sum_{i,j=1}^r \delta X_i(x)S_{ij} V_{j}(v) + X_i(x)\delta S_{ij} V_{j}(v) + X_i(x)S_{ij} \delta V_{j}(v)\\
    &\text{ with } \delta S_{ij} \in\mathbb{R},\delta X_i \in L^2(D_x), \delta V_j\in L^2(D_v) \text{ and } \langle \delta X_i, X_j\rangle_x = 0, \langle \delta V_i, V_j\rangle_{v} = 0 \Bigg\}.
\end{align*}
After establishing the low-rank manifold and its associated tangent space, our objective is to find $f(t,\cdot,\cdot)\in\mathcal{M}_r$ that satisfies two conditions: (1) $\partial_t f(t,\cdot,\cdot) \in\mathcal{T}_f\mathcal{M}_r$, and (2) $\Vert \partial_t f(t,\cdot,\cdot) - F(f(t,\cdot,\cdot))\Vert_{L^2(D_x\times D_v)}$ is minimized, given a specific right-hand side $F(f)$. Defining the projection onto the span of an orthonormal basis $B$ as $P_B$, and the projection onto the corresponding orthogonal space as $P_B^{\perp}$, the problem can be rewritten as \cite{koch2007dynamical} 
\begin{subequations}\label{eq:ill-cond-dlra}
    \begin{align}
        \dot S =& \,\langle F(f(t,\cdot,\cdot)), X(t,\cdot)V(t,\cdot)^{\top}\rangle_{xv}\\ 
        \partial_t X(t,x) =& \,P_X^{\perp} \langle F(f(t,\cdot,\cdot)), V(t,\cdot)^{\top}\rangle_v S^{-1}(t)\\
        \partial_t V(t,v) =& \,P_V^{\perp} \langle F(f(t,\cdot,\cdot)), X(t,\cdot)^{\top}\rangle_x S^{-\top}(t)\, .
    \end{align}
\end{subequations}

\subsection{The {augmented BUG} integrator}\label{sec:rBUG}
Standard time integrators to evolve the system \eqref{eq:ill-cond-dlra} in time require inverting the coefficient matrix $S$, which when $S$ is ill-conditioned leads to prohibitively small time step sizes. A widely used time integrator to robustly solve this system is the {augmented} BUG integrator \cite{ceruti2022rank} (also sometimes known as the rank-adaptive BUG or simply the BUG integrator). Instead of updating the factors $X$ and $V$ in time, it updates $K(t,x) := X(t,x)^{\top}S(t)$ and $L(t,v) := S(t)V(t,v)$ in parallel from a given time $t_0$ to $t_1 = t_0 + \Delta t$. It then retrieves the time-updated basis $X(t_1,x)$ and $V(t_1,x)$ through a Gram-Schmidt process. To enable conservation properties and allow for rank-adaptivity, the integrator augments the time-updated basis with $X(t_0,x)$ and $V(t_0,x)$ and then performs a Galerkin step to update the coefficient matrix $S$. Lastly, the method performs a truncation step according to a user-determined tolerance parameter $\vartheta$ to identify the rank at time $t_1$. In detail, starting with a rank $r$ and the factorized solution at time $t_0$ given by $X^0(x)$, $V^0(v)$, and $S^0$, the algorithm reads as follows:
\begin{enumerate}
\item \textbf{K-Step}: Update and augment the basis $X^0(x)\in\mathbb{R}^r$ to $\widehat X^1(x)\in\mathbb{R}^{2r}$ by solving
\begin{align*}
\partial_t K(t,x) = \left\langle F\left(K(t,x)^{\top} V^0\right), V^0 \right\rangle_v, \quad K(t_0,x) = X^0(x)^{\top} S^0.
\end{align*}
Compute $\widehat X^1(x)\in\mathbb{R}^{2r}$ such that $[X^0, K(t_1,x)]= \widehat X^1(x)^{\top} R$ and store $\widehat M = \langle \widehat X^1, X^{0,\top} \rangle_x\in\mathbb{R}^{2r\times r}$.
\item \textbf{L-Step}: Update and augment the basis $V^0(v)\in\mathbb{R}^r$ to $\widehat V^1(v)\in\mathbb{R}^{2r}$ by solving
\begin{align*}
\partial_t L(t,v) = \left\langle F\left(X^{0,\top} L(t,v)\right), X^0 \right\rangle_x, \quad L(t_0,v) = V^0(v)^{\top} S^{0,\top}.
\end{align*}
Compute $\widehat V^1(v)\in\mathbb{R}^{2r}$ such that $[V^0, L(t_1,v)]= \widehat V^1(v)^{\top} \widetilde R$ and store $\widehat N = \langle \widehat V^1, V^{0,\top} \rangle_v\in\mathbb{R}^{2r\times r}$.
\item \textbf{S-step}: Update $S^0 \in\mathbb{R}^{r\times r}$ to $\widehat S^1\in\mathbb{R}^{2r\times 2r}$ by solving
\begin{align*}
\dot{\widehat S}(t) = \left\langle F\left(\widehat X^{1,\top} \widehat S(t) \widehat V^1\right), \widehat X^1 \widehat V^{1,\top} \right\rangle_{xv}, \quad \widehat S(t_0) = \widehat M S^0 \widehat N^{\top}.
\end{align*}
\item \textbf{Truncation}: Compute the singular value decomposition $\widehat S(t_1) = \widehat P \widehat \Sigma \widehat Q^\top$ where $\Sigma = \text{diag}(\sigma_j)$. For a user-determined tolerance parameter $\vartheta$, choose the new rank $r_1 \leq 2r$ as the minimal $r_1$ such that 
\begin{align*}
\left(\sum_{j=r_1+1}^{2r} \sigma_j^2\right)^{1/2} \leq \vartheta
\end{align*}
holds. Pick $S^1 = \text{diag}(\sigma_1,\cdots,\sigma_{r_1})$ and let $P$ and $Q$ contain the first $r_1$ columns of $\widehat P$ and $\widehat Q$, respectively. Then choose $X^1(x) = \widehat X^{1}(x)^{\top}P\in\mathbb{R}^{r_1}$ and $V^1(v) = \widehat V^{1}(v)^{\top}Q\in\mathbb{R}^{r_1}$.
\end{enumerate}
Then, the factorized solution at time $t_1$ with rank $r_1$ is given as $X^1(x)$, $V^1(v)$, and $S^1$. By successively using the above algorithm, the solution can then be evolved to a given time $t_n$. 

The {augmented} BUG integrator possesses several beneficial properties compared to the fixed-rank BUG integrator of \cite{ceruti2022unconventional}. Perhaps most importantly, besides the ability to achieve rank-adaptivity, the {augmented} BUG integrator fulfills an important conservation property. Unlike the fixed-rank BUG integrator, the rank-adaptive version does not introduce an error when constructing $\widehat S(t_0)$ in the $S$ step equation \cite[Lemma~1]{ceruti2022rank}: Let $f_0 = X^{0,\top} S^0 V^0$ and $\widehat f_0 = \widehat X^{1,\top} \widehat S(t_0) \widehat V^1$. Then, $f_0 = \widehat f_0$, which is due to the fact that the old basis at time $t_0$ is included in the augmented basis used in the Galerkin step. It has been shown in \cite{ceruti2022rank} that this property enables norm conservation as well as energy conservation for Schr\"odinger and Hamiltonian systems. Invariants of kinetic systems such as mass are however not conserved by the augmented BUG integrator. We will show in this paper that the augmentation step guarantees the conservation of invariants up to a user-determined tolerance parameter.

Recent developments of BUG integrators have taken \cite{ceruti2022rank} as a starting point, i.e., these integrators do not introduce errors when constructing $\widehat S(t_0)$. Reviewing all recent modifications to the BUG integrator in this work is not possible, and we refer the reader to \cite{ceruti2023parallel} for the parallel BUG integrator, to \cite{kusch2024second} for its extension to second order, and to \cite{ceruti2024robust} for the extension of the augmented BUG integrator to second order.

\section{Recap: Conservative BUG integrators for DLRA}\label{sec:cBUG}

In the following we focus on partial differential equations that admit a conservation law for a given number of conserved basis functions $U(v)\in\mathbb{R}^m$. The goal is to modify existing BUG integrators to preserve the full-rank dynamics of $\Phi(t,x)= \int_{D_v} f(t,x,v) U(v)\, dv$. That is, the modified BUG integrator should preserve the local conservation law $\partial_t \Phi(t,x) = \langle F(f), U \rangle_v$. Currently, two distinct approaches exist to ensure conservation, namely using modified conservative evolution equations or a micro-macro decompositions. We will summarize these two approaches in the following.

\subsection{The {continuously conservative BUG} integrator of \cite{einkemmer2023robust}}\label{sec:cBUG}
The conservative integrator of \cite{einkemmer2023robust} is constructed to preserve local conservation laws, specifically for the Vlasov equation, however, the strategy extends to other kinetic problems. One core building block is to evolve $V(t,v)$ such that basis functions $U(v)\in\mathbb{R}^m$ which are important for conservation are spanned by the basis. Following \cite{einkemmer2021mass}, the DLRA evolution equations \eqref{eq:ill-cond-dlra} are modified to preserve the basis directions $U(v)$. Denoting $W_p(t,v) = V_p(t,v)$ for $m < p\leq r$, the DLRA evolution equations with $f_r(t,x,v) = X(t,x)^{\top} S(t) V(t,v)$ change to
\begin{subequations}\label{eq:modifiedDLRA}
    \begin{alignat}{2}
    \sum_{i = 1}^r \dot X_i S_{ik} =& \langle V_k, F(f_r)\rangle_v - \sum_{i = 1}^r X_i \dot S_{ik}, \qquad &&1 \leq k \leq r\\
    \sum_{i = 1}^r\sum_{p = m+1}^r S_{iq} S_{ip} \dot{W}_p =& \sum_{i=1}^rS_{iq} \langle S_{iq}\langle X_i, F(f_r)\rangle_x - \sum_{i,l = 1}^r S_{iq}\dot S_{il} V_l, \qquad &&m+1 \leq q \leq r\\
    \dot S_{kl} =& \langle X_k V_l, F(f_r) \rangle_{xv}, \qquad &&1\leq k,l\leq r.
\end{alignat}
\end{subequations}

Here, we use index notation to precisely define the ranges of sums and leave out dependencies on the phase space for better readability. While this system has been solved through a non-robust time integration scheme in \cite{einkemmer2021mass}, the BUG integrator has been employed in \cite{einkemmer2023robust} to obtain stability independent of the curvature of the low-rank manifold. To remove the projection error of the $S$ step in the fixed-rank integrator, the augmentation strategy of \cite[Lemma~1]{ceruti2022rank} is used. That is, the basis at the previous time step $t_0$ is included in the Galerkin step to remove the projection error of the initial condition for the Galerkin step. 

{Moreover, the authors add the basis $\nabla X^0$ in order to be able to represent the local conservation law $\partial_t \rho + \nabla \cdot j = 0$ exactly. That $\nabla \cdot j$ lies in the approximation space then also implies that we can represent the local conservation laws for momentum and energy exactly (on the continuous level). However, what the authors in \cite{einkemmer2023robust} did not realize is that there is no need to explicitly add $\nabla X^0$ to the basis. To see this it is sufficient to consider the local conservation law for mass. In the K-step}
\[  \nabla \cdot j = -\partial_t \rho = \int \partial_t f \,dv  = \int F(f) \,dv = \langle F(f), 1\rangle_v \]
{Thus, $\nabla \cdot j$ is already ensured to lie in $\text{span}\{X^1\}$ as it is ensured that $1$ lies in $\text{span}\{V^0\}$. This has two important consequences. First, adding the additional basis functions results in a larger space and thus increased computational cost. Second, in some situations, the fact that we add redundant basis functions to the orthogonalization procedure can result in increased round-off error. We thus, will denote the variant of the algorithm by \textit{continously conservative integrator} which does not add $\nabla X^0$ to the basis.} 

The continuously conservative BUG integrator differs in two points from the augmented BUG integrator \cite{ceruti2022rank}. First, it employs a conservative truncation. Second, the continuously conservative BUG integrator solves system \eqref{eq:modifiedDLRA} instead of \eqref{eq:ill-cond-dlra}.

\subsection{The macro-micro decomposition of \cite{koellermeier2024macro}}
The conservative integrator of \cite{koellermeier2024macro} is constructed to preserve local conservation laws, specifically for kinetic shallow-water moment models, however, the strategy extends to other kinetic problems. The main idea is to define a solution representation $f(t,x,v) \approx \varphi(t,x,v) + \psi(t,x,v)$, where $\varphi$ is the conserved solution which determines the evolution of $\Phi$ and $\psi$ is the remainder such that $\langle \psi, U \rangle_v = 0$. Representations of $\varphi$, for example, include $\varphi(t,x,v) = \Phi(t,x)^{\top}U(v)$. A conventional method is then used to evolve the conserved quantities $\Phi(t,x)$, preserving the local conservation law
\begin{align*}
    \partial_t \Phi(t,x) = \langle F(\varphi + \psi), U \rangle_v\,.
\end{align*}
The remainder is then evolved with a DLRA integrator ensuring that $\langle \psi, U \rangle_v = 0$ at all times. A similar idea has been used for the Vlasov-Poisson systems in \cite{coughlin2024robust} (using a nodal velocity discretization and the projector--splitting integrator) and in \cite{frank2025asymptotic} for thermal radiative transfer. While micro-macro approaches help ensure structure--preservation, they commonly require more significant modifications to an existing code base and, in the case of nodal representations, require additional orthonormalizations. Unlike the approach of \cite{einkemmer2023robust} they can, however, be combined with more general time integration schemes and the extension to higher--order BUG integrators is straightforward.

\section{Construction of high order conservative BUG integrators}\label{sec:local-cons-rbug}
In the following, we demonstrate a strategy to ensure local conservation properties for general BUG integrators when using explicit Runge-Kutta methods to solve the coefficient update. As the name suggests, most basis-update \& Galerkin integrators consist of a step that updates the basis functions, followed by a step to update the coefficient matrix $S$, usually called the $S$-step. In that case, the dynamics are solely determined in the $S$-step which is, therefore, the key component to studying local conservation. We note that parallel BUG integrators evolve the dynamics not only in a sequential coefficient update, an observation that we will discuss in the next section and which will actually help in the construction of conservative methods.

Given a vector of conserved basis functions $U(v)\in\mathbb{R}^m$, our goal is to construct a method that fulfills local conservation laws of the form
\begin{align}\label{eq:localConsBUG}
    \partial_t \phi(t,x) =&\, \left\langle F\left(f(t,x,v)\right), U(v) \right\rangle_{v}\,,
\end{align}
where $\phi(t,x) := \langle f(t,x,v), U(v) \rangle_{v}$. Conventional BUG integrators do not satisfy \eqref{eq:localConsBUG} for two distinct reasons. First, a truncation step is required between time updates. This truncation is commonly done by SVD, which violates conservation. To remedy this, the conservative truncation step in \cite{einkemmer2023robust,guo2022conservative} can be used. Such a truncation step ensures that if the vector of $m$ conserved basis functions $U(v)$ lies in the augmented basis, i.e., $P_{\widehat V^1} U(v) = U(v)$, then denoting the basis after truncation as $V^1$, we have $P_{V^1} U(v) = U(v)$. We note that when $U(v)$ lies in the velocity basis at time $t=0$ and the augmented basis incorporates the basis at the old time step (which is the case for all commonly used BUG integrators), then $P_{V^n} U(v) = P_{\widehat V^n} U(v) = U(v)$ for all time steps $n$. Second, the spatial basis might not capture the correct dynamics of the local conservation law \eqref{eq:localConsBUG}. To illustrate this, let us compare the $S$-step to our local conservation law. Let us denote $\widehat X = \widehat X^{n+1}(x)\in \mathbb{R}^{\widehat r}$, $\widehat V = \widehat V^{n+1}(v)\in \mathbb{R}^{\widehat r}$ as the basis of a general BUG integrator, i.e., we omit dependence on the phase space and the superscript $n+1$ in the following discussion for ease of notation. With $Y_S(t,x,v):=\widehat X(x)^{\top}\widehat S(t)\widehat V(v)$, the $S$-step can be written as
\begin{align*}
    \partial_t Y_S(t,x,v) = P_{\widehat V}P_{\widehat X}F(Y_S(t,x,v))\,. 
\end{align*}
To arrive at a local conservation law, we multiply the above equation with our vector of conserved basis functions $U(v)$ and integrate over velocity. With $\phi_S(t,x):= \langle \widehat X(x)^{\top} \widehat S(t)\widehat V(v), U(v)\rangle_v$, this yields
\begin{align}
    \partial_t \phi_S(t,x) =\,& \langle P_{\widehat V}P_{\widehat X}F(Y_S(t,x,v)), U(v)\rangle_v\nonumber \\
    =\,& \langle P_{\widehat X}F(Y_S(t,x,v)), P_{\widehat V} U(v)\rangle_v\nonumber  \\
    =\,& P_{\widehat X}\langle F(Y_S(t,x,v)), U(v)\rangle_v\,, \label{eq:phiScont}
\end{align}
where we have used that $P_{\widehat V}$ is self-adjoint and $P_{\widehat V} U(v) = U(v)$ when using a conservative truncation as discussed above. This expression is similar to the exact local conservation law \eqref{eq:localConsBUG}, however, it incorporates the projection onto $\widehat X$, which still prevents us from obtaining the conservation law. In particular, it prevents us from integrating equation \eqref{eq:phiScont} in space to obtain global conservation. To eliminate the projection, we have to ensure that $\langle F(Y_S(t,x,v)), U(v)\rangle_v$ is spanned by the augmented spatial basis at all times $t$, which is in general not possible. However, we note that this conservation law only needs to be fulfilled at the numerical level after a time discretization is performed. Hence, in the case of a general explicit $s$-stage Runge-Kutta method (given by coefficients $a_{ij}$ and $b_i$) we have for the $S$-step
\begin{subequations}\label{eq:S-rk-pre-all}
\begin{align}
    \widehat S^{n+1} = \widehat S^{n} + \Delta t \sum_{i=1}^{s} b_i F_S^{(i)}, \qquad \text{with} \qquad  \widehat S^n := \langle \widehat X, X^{n,\top} \rangle_{x} S^n \langle V^n, \widehat V^{\top} \rangle_{v}, \label{eq:S-rk-pre}
\end{align}
where the stages $F_S^{(i)}$ are computed as
\begin{align}
F_S^{(i)} = \left\langle F^{(i)}, \widehat X\widehat V^{\top} \right\rangle_{x,v}, \quad \text{with }F^{(i)} := F \left( \widehat X^{\top} \widehat S^n \widehat V + \Delta t \sum_{j=1}^{i-1} a_{ij} \widehat X^{\top}F_S^{(j)} \widehat V \right).
\end{align}
\end{subequations}
Our main goal is to enforce a discrete counterpart of the local conservation law \eqref{eq:localConsBUG} which takes the form
\begin{align}\label{eq:localConsBUGdiscrete}
    \phi^{n+1}(x) = \phi^{n}(x) + \Delta t \sum_{i=1}^{s} b_i \langle F^{(i)}, U(v)\rangle_v \,.
\end{align}
For $\phi^n_S(x):= \langle \widehat X(x)^{\top} \widehat S^n \widehat V(v), U(v)\rangle_v$ the discrete counterpart of \eqref{eq:phiScont} however reads
\begin{align}\label{eq:localConsBUGdiscreteP}
    \phi_S^{n+1}(x) = \phi_S^{n}(x) + \Delta t P_{\widehat X} \sum_{i=1}^{s} b_i \langle F^{(i)}, U(v)\rangle_v \,.
\end{align}
This is similar to the desired conservation law \eqref{eq:localConsBUGdiscrete}, however, just as in the time-continuous case, the equation incorporates a projection term $P_{\widehat X}$. Therefore, to enforce conservation, we need to ensure that $\sum_{i=1}^{s} b_i \langle F^{(i)}, U(v)\rangle_v$ is spanned by the spatial basis. We can enforce this through an additional basis augmentation step. That is, if we denote an orthonormalization step (e.g. modified Gram-Schmidt) by $\mathrm{ortho}$, then we define the new basis
\begin{align*}
    \bar X = [\widehat X, \bar X^{\star}] := \mathrm{ortho}\left(\widehat X, \sum_{i=1}^{s} b_i \langle F^{(i)}, U(v)\rangle_v\right) \in\mathbb{R}^{\widehat r + m}.
\end{align*}
Using this new basis, the Runge-Kutta update of the $S$-step becomes 
\begin{align}\label{eq:S-rk}
    \bar S^{n+1} = \bar S^{n} + \Delta t \sum_{i=1}^{s} b_i \bar F_S^{(i)} \qquad \text{with} \qquad \bar S^n := \langle \bar X, \widehat X^{n,\top} \rangle_{x} \widehat S^n,
\end{align}
where the stages $\bar F_S^{(i)}$ are computed as
\begin{align*}
\bar F_S^{(i)} = \left\langle F^{(i)}, \bar X\widehat V^{\top} \right\rangle_{x,v}, \quad \text{with }F^{(i)} := F \left( \widehat X^{\top} \widehat S^n \widehat V + \Delta t \sum_{j=1}^{i-1} a_{ij} \widehat X^{\top}F_S^{(j)} \widehat V \right).
\end{align*}
By construction of the basis $\bar X$, when redefining $\phi^n_S(x):= \langle \bar X(x)^{\top} \bar S^n \widehat V(v), U(v)\rangle_v$ we get
\begin{align*}
    \phi_S^{n+1}(x) =\,& \phi_S^{n}(x) + \Delta t P_{\bar X} \sum_{i=1}^{s} b_i \langle F^{(i)}, U(v)\rangle_v \,\\
    =\,& \phi_S^{n}(x) + \Delta t \sum_{i=1}^{s} b_i \langle F^{(i)}, U(v)\rangle_v\,.
\end{align*}
Hence, we fulfill the discrete local conservation law \eqref{eq:localConsBUGdiscrete}. This then also trivially implies global conservation as $\langle F^{(i)}, U(v)\rangle_{x,v} = 0$ (a property of the continuous model) and thus
\[ \int \phi_S^{n+1} \,dx = \int \phi_S^{n} \,dx.\]

Let us note that we have constructed prohibitively large objects $F^{(i)}$ in our discussion, which we do not want to store and compute. Instead, it is sufficient to only compute the projections $F_S^{(i)} = \langle F^{(i)}, \widehat X(x)\widehat V(v)^{\top}\rangle_{xv}$ and $\bar F_{\star}^{(i)} := \langle F^{(i)}, \bar X^{\star}(x)\widehat V(v)^{\top}\rangle_{xv}$, where $\bar F_S^{(i)} = [F_S^{(i)}, \bar F_{\star}^{(i)}]^{\top}$ instead. To reduce the computation of $\bar F_X^{(i)}$, we can instead compute and store $\bar F_X^{(i)} := \langle F^{(i)}, \widehat V(v)^{\top}\rangle_{v}$ and determine $\bar F_{\star}^{(i)} = \langle \bar F_X^{(i)}, \bar X^{\star}(x)\rangle_{x}$. Let us, in the following, write down the integrator that followed from our discussion.

\begin{enumerate}
\item \textbf{K-Step}: Determine an augmented basis $\widehat X^1(x)\in\mathbb{R}^{\widehat r}$ which spans the old basis $X^{0}(x)$ with the $K$-step of any BUG integrator (except for the fixed-rank BUG integrator).
\\[2mm]
\item \textbf{L-Step}: Determine an augmented basis $\widehat V^1(v)\in\mathbb{R}^{\widehat r}$ which spans the old basis $V^{0}(v)$ with the $L$-step of any BUG integrator (except for the fixed-rank BUG integrator).\\[2mm]
\item \textbf{S-step}: Update $S^0 \in\mathbb{R}^{r\times r}$ to $\widehat S^1\in\mathbb{R}^{(\widehat r + m)\times \widehat r}$. Proceed as follows: 
\begin{enumerate}
    \item For $i = 1,\cdots, s$ compute and store the stages $F_S^{(i)}$ and $\bar F_X^{(i)}$ as
        \begin{align*}
            F_S^{(i)} = \left\langle \bar F_X^{(i)}, \widehat X \right\rangle_{x}, \quad \text{with } \bar F_X^{(i)} = \langle F^{(i)}, \widehat V^{\top}\rangle_v .
        \end{align*}
    \item Compute and store the augmented basis 
        \begin{align*}
            \bar X = [\widehat X, \bar X^{\star}] := \mathrm{ortho}\left(\widehat X, \sum_{i=1}^{s} b_i \bar{F}_X^{(i)} \langle \widehat{V}, U^{\top}\rangle_v\right) \in\mathbb{R}^{\widehat r + m}\,.
        \end{align*}
    \item Compute and store the augmented stages 
        \begin{align*}
            \bar F_S^{(i)} = [F_S^{(i),\top}, \bar F_{\star}^{(i),\top}]^{\top} \in\mathbb{R}^{(\widehat r + m)\times \widehat r }, \quad \text{with } \bar F_{\star}^{(i)} = \langle \bar F_X^{(i)}, \bar X^{\star}(x)\rangle_{x}\,.
        \end{align*}
    \item Update coefficients with
        \begin{align}\label{eq:S-rk}
            \bar S^{n+1} = \bar S^{n} + \Delta t \sum_{i=1}^{s} b_i \bar F_S^{(i)}\,.
        \end{align}
\end{enumerate}
\item \textbf{Truncation}: Perform a conservative truncation of $\bar Y^{1} = \bar X^{1, \top}\bar S^{1} \widehat{V}^{1}$ to a rank $r_{1}$ solution $Y^{1} = X^{1, \top} S^{1} V^{1}$ given a user-determined tolerance parameter $\vartheta$. 
\end{enumerate}

The strategy proposed in this section requires an additional QR factorization to evaluate $\bar X_{\star}$, which, however, in most situations is unlikely to be a bottleneck. We also note that the augmentation is specific to the Runge--Kutta method used. Thus, changing the time integrator also requires a change in $\bar{X}$ to maintain conservation. The outlined procedure can, however, be applied to any (commonly used) BUG-type integrator, with the exeception of the original fixed-rank BUG integrator of \cite{ceruti2022unconventional}. 

\section{Examples}\label{sec:examples}
In the following, we discuss two examples where the scheme is conservative without requiring additional augmentations. We recall that the two key ingredients to local conservation are 1) the preservation of conserved basis functions $U(v)$ which is directly fulfilled for BUG integrators with conservative truncation, and 2) ensuring that the spatial basis fulfills (see equation \eqref{eq:localConsBUGdiscreteP}) 
\[ P_{\widehat X}  \sum_{i=1}^s b_i \langle F^{(i)}, U(v)\rangle_v = \sum_{i=1}^s b_i \langle F^{(i)}, U(v)\rangle_v. \]
In the following, we show that these conditions are automatically fulfilled when combining the augmented BUG and the midpoint BUG integrators with specific time integration methods in the substeps. Thus, in this case, the time integrator used to solve the substeps is compatible with the robust DLRA integrator.

\subsection{Augmented BUG integrator of \cite{ceruti2022rank}}
First, we focus on the augmented BUG integrator in combination with an explicit Euler time updates in the $K$ and $S$-steps. In this case, the Runge--Kutta stage in \eqref{eq:S-rk-pre-all} becomes $F^{(1)} = F(Y^n(x,\cdot))$ and the corresponding $K$-step from time $t_0$ to $t_1$ reads
\begin{align*}
    K^{n+1}(x) = K^n(x) + \Delta t \langle F(Y^n(x,\cdot)), V^n \rangle_v = K^n(x) + \Delta t \langle F^{(1)}, V^n \rangle_v\,.
\end{align*}
We can then show that $\langle F^{(1)},U(v)\rangle_v$ lies in the range of our basis $\widehat X^{n+1}$. Note that
\begin{align*}
   \langle F^{(1)},U\rangle_v = \sum_i\left\langle F^{(1)},\langle U,V_i^n\rangle_v V_i^n\right\rangle_v  = \frac1{\Delta t}\sum_i (K_i^{n+1}(x) - K_i^{n}(x))\langle U, V_i^n\rangle_v.
\end{align*}
Since in the augmentation step $[X^{n}(x), K^{n+1}(x)]= \widehat X^{n+1}(x)^{\top} R$ the term $K_i^{n+1}(x) - K_i^{n}(x)$ is spanned exactly by the augmented basis $\widehat X^{n+1}(x)$, we know that 
\begin{align*}
    P_{\widehat X^{n+1}}\langle F^{(1)},U\rangle_v = \langle F^{(1)},U\rangle_v.
\end{align*}
Therefore, \eqref{eq:localConsBUGdiscreteP} when using a forward Euler time discretization becomes with $b_1 = 1$
\begin{align*}
    \phi_S^{n+1}(x) = \phi_S^{n}(x) + \Delta t P_{\widehat X} b_1 \langle F^{(1)}, U(v)\rangle_v = \phi_S^{n}(x) + \Delta t \langle F^{(1)}, U(v)\rangle_v \,.
\end{align*}
The augmented BUG integrator is hence conservative when using a forward Euler time discretization in the $K$ and $S$-steps combined with a conservative truncation. An additional augmentation is not required.
\subsection{Midpoint BUG integrator of \cite{ceruti2024robust}}
We investigate the midpoint BUG integrator in combination with a midpoint Runge--Kutta scheme in the $S$-step. In this case, the augmented spatial basis at which the coefficient update is computed takes the form
\begin{align*}
    \widehat X = \text{ortho}(X^n, K(t_{n+1}), \langle F(\widehat Y^{n+1/2}), \widehat V^{n+1/2}\rangle_v)\,,
\end{align*}
where the rank $2r$ solution $\widehat Y^{n+1/2}$ and velocity basis $\widehat V^{n+1/2}$ can be computed with the augmented BUG integrator (without a truncation step). To ensure conservation while preserving second-order accuracy (again, the latter formally requires a Lipschitz continuous and bounded right-hand side), we choose a time discretization of the $S$-step with Runge-Kutta stages $F^{(1)}$ (determined by the time integration method to compute $\widehat Y^{n+1/2}$) and $F^{(2)} = F(\widehat Y^{n+1/2})$ as well as $b_1=0$ and $b_2 = 1$, see \eqref{eq:S-rk-pre-all}. Then,
\begin{align*}
   \langle F^{(2)},U\rangle_v = \sum_i\left\langle F^{(2)},\langle U,\widehat V_i^{n+1/2}\rangle_v \widehat V_i^{n+1/2}\right\rangle_v = \sum_i\left\langle F^{(2)}, \widehat V_i^{n+1/2}\right\rangle_v \langle U,\widehat V_i^{n+1/2}\rangle_v\,.
\end{align*}
Due to the choice of $\widehat X$ and $F^{(2)} = F(\widehat Y^{n+1/2})$, we have
\begin{align*}
   P_{\widehat X}\langle F^{(2)},U\rangle_v =\,& \sum_i P_{\widehat X}\left\langle F^{(2)}, \widehat V_i^{n+1/2}\right\rangle_v  \langle U,\widehat V_i^{n+1/2}\rangle_v\\
   =\,& \sum_i \left\langle F^{(2)}, \widehat V_i^{n+1/2}\right\rangle_v  \langle U,\widehat V_i^{n+1/2}\rangle_v = \langle F^{(2)},U\rangle_v\,.\\
\end{align*}
Therefore, the midpoint BUG integrator with a midpoint Runge--Kutta time discretization in the $S$-step which uses $F^{(2)} = F(\widehat Y^{n+1/2})$ combined with a conservative truncation, is locally conservative. An additional augmentation is not required.

\section{A conservative second-order parallel BUG integrator}\label{sec:cons-parallel}

Parallel BUG integrators do not require the computation of a sequential $S$-step. Instead, the $S$-step is computed as part of the basis update and the updated coefficient matrix is generated in a computationally inexpensive coefficient augmentation step. This helps us to construct parallel conservative BUG integrators. Let us directly state the integrator and then show that it is locally conservative by following the discussion of the previous section. 
\begin{enumerate}
\item Construct augmented basis matrices $\widehat X^0 = \text{ortho}([X^0, \langle F(Y_0), V^0\rangle_v])$ and $\widehat V^0 = \text{ortho}([ V^0, \langle F(t_0, Y_0), X^0\rangle_x])$.
\item Determine the augmented basis matrices $\widehat{U}_1\in \mathbb{R}^{m\times 4r}$ and $\widehat{V}_1\in \mathbb{R}^{n\times 4r}$(in parallel):
\\[2mm]
\textbf{K-step}: Integrate from $t=t_0$ to $t_1$ the $m \times 2r$ matrix differential equation
\begin{align*}
\dot{K}(t) =\,& \langle F(t, K(t)^{\top}\widehat  V^0), \widehat  V^0\rangle\, , \quad  K(t_0) = X^{0,\top}  S_0  \langle V^0, \widehat{ V}_0^{\top}\rangle_v\, .
\end{align*}
Compute $\widehat X^1(x) = [\widehat X^0, \widetilde X^2]\in\mathbb{R}^{4r}$ such that $[X^0, K(t_1,x)]= \widehat X^{1,\top} R_1$.
\\[2mm]
\textbf{L-step}: Integrate from $t=t_0$ to $t_1$ the $n \times 2r$ matrix differential equation
\begin{align*}
\dot{L}(t) =\,& \langle F(t,\widehat X^{0,\top}  L(t)), \widehat  X^0\rangle_x\, , \quad  L(t_0) =  V^{0,\top}  S_0^{\top} \langle X^0,\widehat{X}_0^{\top}\rangle_x\, .
\end{align*}
Compute $\widehat V^1(v) = [\widehat V^0, \widetilde V^2]\in\mathbb{R}^{4r}$ such that $[V^0, L(t_1,v)]= \widehat V^{1,\top} R_2$.
\item \textbf{Augment}: Set up the augmented coefficient matrix $\widehat{ S}_1 \in\mathbb{R}^{4r\times 4r}$ as
\begin{align*}
\widehat{ S}_1 = \begin{pmatrix}
\langle \widehat X^{0,\top},  K(t_1,x)\rangle_x &  \langle L(t_1,v)^{\top}, \widetilde  V^2\rangle_v\\
\langle \widetilde X^{2,\top},  K(t_1,x)\rangle_x & \bf0
\end{pmatrix}\, .
\end{align*}
\item \textbf{Truncate}: Perform a conservative truncation of $\widehat Y^{1} = \widehat X^{1, \top}\widehat S^{1} \widehat{V}^{1}$ to a rank $r_{1}$ solution $Y^{1} = X^{1, \top} S^{1} V^{1}$ given a user-determined tolerance parameter $\vartheta$. 
\end{enumerate}

The integrator follows the second-order parallel integrator of \cite{kusch2024second}, but replaces the $S$-step by the $K$-update. That is, the $S$-step is not computed and the corresponding block in the augmentation step is replaced by $\langle \widehat X^{0,\top},  K(t_1,x)\rangle_x$. This modification preserves the second-order error bound of \cite{kusch2024second} under the same assumptions. This is the consequence of the $K$-step carrying second order information for the bases $\widehat X^{0}$ and $\widehat V^{0}$, which follows from a trivial modification of the arguments in \cite[Lemma~4.4]{kusch2024second}.
It is easily verified that the stated integrator is locally conservative. Again, due to the conservative truncation, we have $P_{\widehat V^1}U = P_{\widehat V^0}U = P_{V^0}U = U$ and $\langle \widetilde V^2, U\rangle = 0$. Then, for $\phi^1_S(x) := \widehat X^{\top}\widehat S\langle \widehat V, U\rangle_v$ we have
\begin{align*}
    \phi^1_S(x) =\,& P_{\widehat X^0}K(t_1,x)^{\top}\langle \widehat V^0, U\rangle_v + P_{\widetilde X^2}K(t_1,x)^{\top}\langle \widehat V^0, U\rangle_v + \langle P_{\widetilde V^2}X^{0,\top}L(t_1,v)\widehat V^0, U\rangle_v \\
    =\,& P_{\widehat X}K(t_1,x)^{\top}\langle \widehat V^0, U\rangle_v\\
    =\,& K(t_1,x)^{\top}\langle \widehat V^0, U\rangle_v\,.
\end{align*}
Plugging in the equation for $K$, and defining $Y_K(t):= K(t)^{\top}\widehat  V^0$ we have
\begin{align*}
    \phi^1_S(x) =\,& \phi^0_S(x) + \int_{t_0}^{t_1} \langle F(t, Y_K(t)), \widehat  V^0\rangle \langle \widehat V^0, U\rangle_v\,dt\\
    =\,& \phi^0_S(x) + \int_{t_0}^{t_1} \langle F(t, Y_K(t)), P_{\widehat  V^0} U\rangle\,dt\\
    =\,& \phi^0_S(x) + \int_{t_0}^{t_1} \langle F(t, Y_K(t)), U\rangle\,dt\,.
\end{align*}
Thus, the dynamics of the invariances are determined by the local conservation law
\begin{equation} \label{eq:cont-cons-YK}
   \partial_t \phi_S(x) = \langle F(t, Y_K(t)), U\rangle_v \,.
\end{equation}

Note that the above calculations show that we obtain a conservation law. Thus, the dynamical low-rank integrator is conservative. We note, however, that we do not exactly reproduce the continuous conservation law due to the presence of $Y_K$ instead of $Y_S$ on the right-hand side. Thus, in general, an error on the order of the low-rank approximation is committed for the flux function.

The above calculations are done for the continuous case. However, the argument carries over to arbitrary Runge--Kutta methods. In that case, we have
\begin{align*}
    \phi_S^1(x) &= K^1(x)^{\top}\langle \widehat V^0, U\rangle_v \\
                &= \phi_S^0(x) + \Delta t \sum_{i=1}^{s} b_i \langle F(Y_K^{(i)}), \widehat V^0\rangle_v  \langle \widehat V^0, U\rangle_v \\
                &= \phi_S^0(x) + \Delta t \sum_{i=1}^{s} b_i \langle F(Y_K^{(i)}), U\rangle_v,
\end{align*}
which is precisely the Runge--Kutta method applied to the continuous conservation law in equation \eqref{eq:cont-cons-YK}. Note that $Y_K^{(i)}$ are the stages of the Runge--Kutta scheme.  We also have global conservation since all Runge--Kutta methods preserve linear invariants.

Let us remark two important details for both algorithms presented so far:
\begin{enumerate}
    \item Note that we assume existence of integrals $\int U(v)\,dv$. When $U(v)$ does not have bounded support, a weighted integral can be defined according to \cite{einkemmer2021mass} to ensure existence.
    \item To use the conservative truncation as described in Section~\ref{sec:cBUG}, it is important to perform the augmentation as described in the above algorithm. That is we augment according to $[V^0,L^1]= \widehat V^1(v)^{\top} \widetilde R$ instead of $[L^1, V^0]= \widehat V^1(v)^{\top} \widetilde R$. In this case, we ensure that the conservative basis $U(v)$ is stored as the first $m$ basis functions in $\widehat V^1(v)$.
\end{enumerate}

\section{Conservation for the linear Radiative Transfer Equation}\label{sec:results_rad}
In this section, we investigate two radiation transport test cases for which we compare the integrators presented in this work. All numerical experiments can be reproduced with the openly available source code\footnote{\url{https://github.com/ScSteffen/Publication_High_Order_Conservative_Integrators.git}}.

In particular, we compare the following conservative methods for the plansesource and linesource test cases:
\begin{itemize}
  \item \textbf{cons. BUG+Euler}: We use the augmented BUG of \cite{ceruti2022rank} that integrates the K,L, and S step with an explicit Euler scheme and is equipped with a conservative truncation. This method does not require any additional augmentation to be conservative.
  \item \textbf{cons. Parallel+Euler}: We use the parallel integrator of \cite{ceruti2023parallel} that integrates the K,L, and S step with an explicit Euler scheme and is equipped with a conservative truncation. This method does not require any additional augmentation to be conservative.
   \item \textbf{cons. 2nd order Parallel+RK4}: We use the second order accurate parallel integrator of \cite{kusch2024second} based on the midpoint rule  that integrates the K,L, and S step with an explicit RK4 scheme and is equipped with a conservative truncation. This method does not require any additional augmentation to be conservative.
  \item \textbf{cons. midp BUG+Euler}: We use the second order integrator of \cite{ceruti2024robust} based on the midpoint rule, that integrates the $K$,$L$, and $S$-step with an explicit Euler scheme and is equipped with a conservative truncation.
  \item \textbf{cons. midp BUG+Heun}: We use the second order integrator of \cite{ceruti2024robust} based on the midpoint rule, that integrates the K, and L step with an explicit Heun scheme, the S step with an explicit RK4 scheme and is equipped with a conservative truncation. This method does not require any additional augmentation to be conservative.
    \item \textbf{cons. midp BUG+RK4}: We use the second order integrator of \cite{ceruti2024robust} based on the midpoint rule, that integrates the K, L and S step with an explicit RK4 scheme and is equipped with a conservative truncation.
\end{itemize}

\begin{table}[t]
\centering
\caption{Overview of integrators compared across test cases. The first group includes the proposed conservative schemes. The second group comprises the non-conservative augmented BUG integrator from \cite{ceruti2022rank}, which serves as a baseline method.}
\label{tab_integrators}
\resizebox{0.5\textwidth}{!}{
\begin{tabular}{@{}l|cc@{}}
\toprule
& \multicolumn{2}{c}{radiation transport} \\
 &1D planesource & 2D linesource  \\
 \midrule
cons. BUG+Euler &\cmark\,\,\, & \cmark\\
cons. Parallel+Euler &\cmark\,\,\, & \xmark\\
cons. 2nd order Parallel+RK4 &\cmark\,  & \cmark\\
cons. midp BUG+Euler &\cmark\,  & \xmark \\
cons. midp BUG+Heun  &\cmark\,  & \xmark\\
cons. midp BUG+RK4   &\cmark\,  & \cmark \\
 \midrule
 non. cons. BUG of \cite{ceruti2022rank} &\cmark\,\,\, &\xmark\\
\bottomrule
\end{tabular}
}
\end{table}

\begin{figure}[t]
    \centering
    \begin{subfigure}{0.32\textwidth}
        \includegraphics[width=\textwidth]{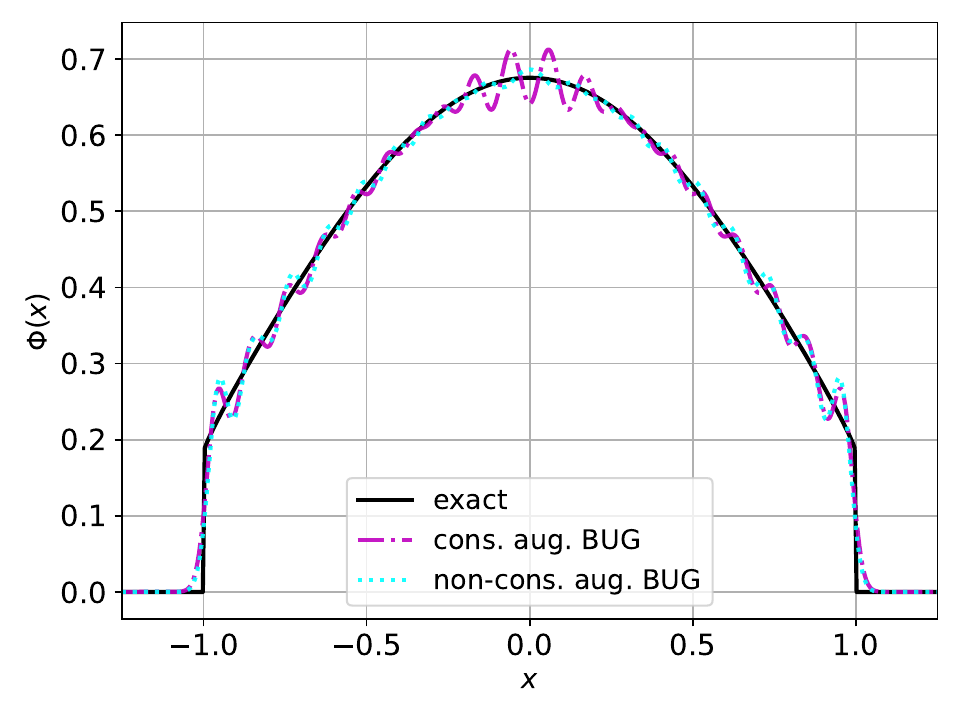}
        \caption{Augmented BUG integrators}
        \label{fig:scalar_flux_dim1_a}
    \end{subfigure}
    \begin{subfigure}{0.32\textwidth}
        \includegraphics[width=\textwidth]{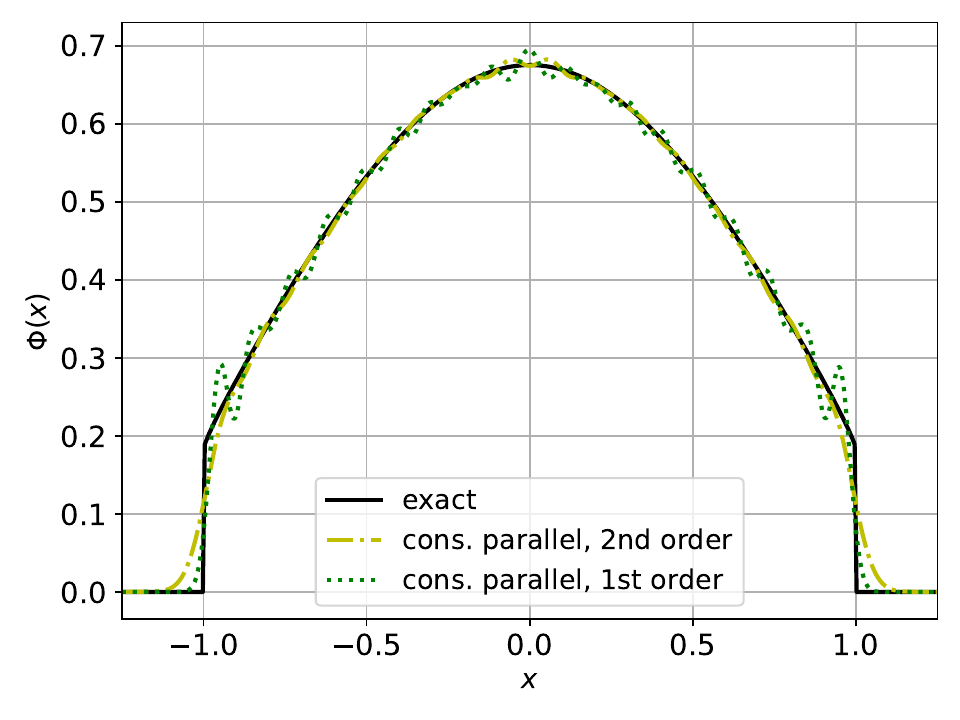}
        \caption{Parallel Integrators}
        \label{fig:scalar_flux_dim1_b}
    \end{subfigure}
    \begin{subfigure}{0.32\textwidth}
        \includegraphics[width=\textwidth]{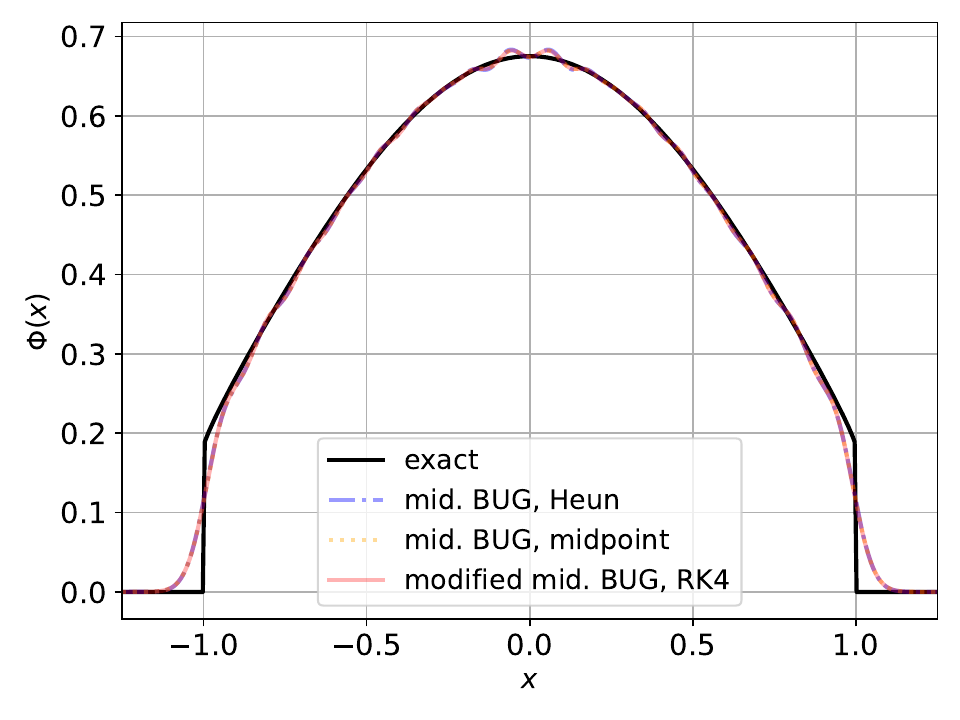}
        \caption{Midpoint BUG integrators}
        \label{fig:scalar_flux_dim1_c}
    \end{subfigure}\\
        \centering
    \begin{subfigure}{0.32\textwidth}
        \includegraphics[width=\textwidth]{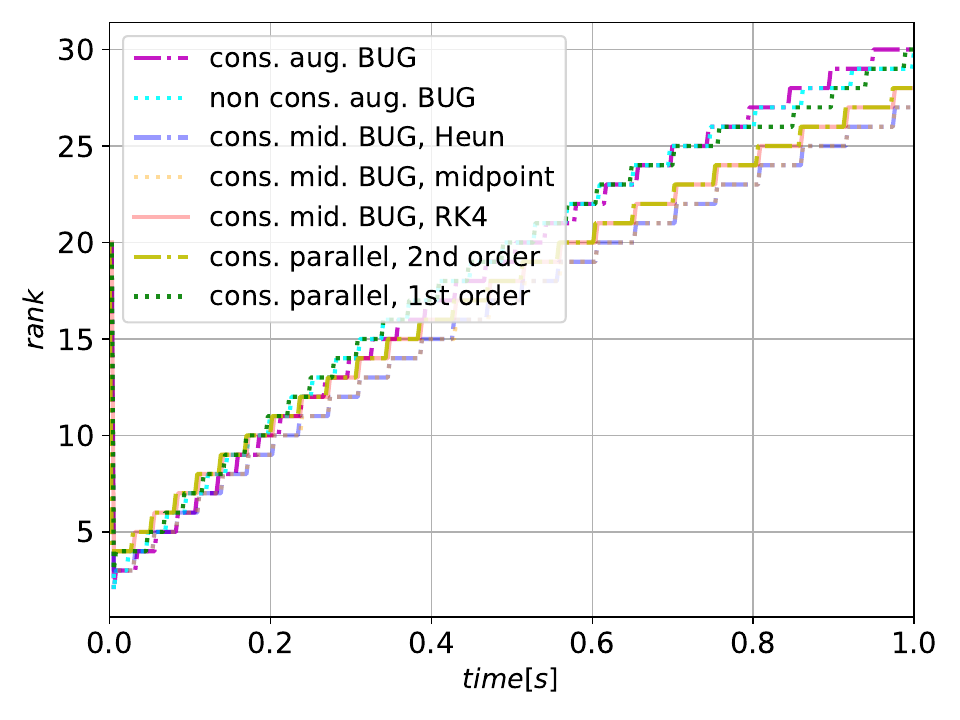}
        \caption{Rank dynamics}
        \label{fig:rank_flux_dim1_a}
    \end{subfigure}
      \begin{subfigure}{0.32\textwidth}
        \includegraphics[width=\textwidth]{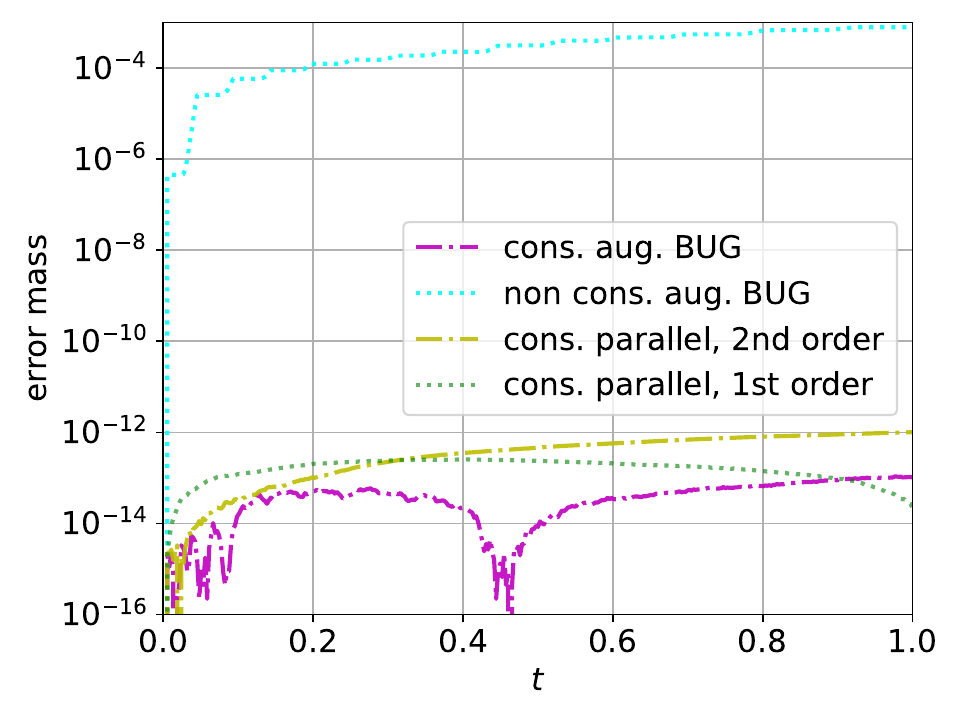}
        \caption{Mass error, BUG and par. integrators}
        \label{fig:mass_flux_dim1_c}
    \end{subfigure}
       \begin{subfigure}{0.32\textwidth}
        \includegraphics[width=\textwidth]{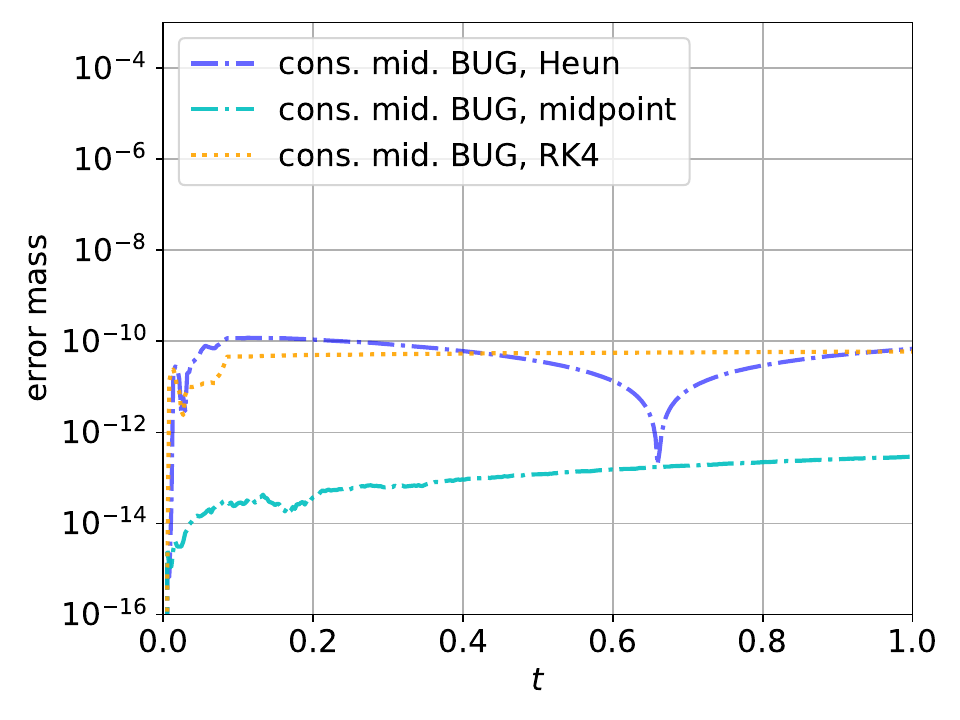}
        \caption{Mass error, midpoint integrators}
        \label{fig:mass_flux_dim1_b}
    \end{subfigure}
    \caption{1D1V Planesource. Comparison of conservative integrators and the non-cons integrator for the line source test case at  $t_{\mathrm{end}}= 1$. All integrators are rank adaptive with threshold $\vartheta=\bar{\vartheta} \Vert \widehat S^1 \Vert$ with $\bar\vartheta = 5\cdot 10^{-2}$. All methods exhibit similar solution behavior (top row) and rank dynamics (bottom left). Midpoint integrators and the second order parallel integrator are less oscillatory, but smear out slightly more. The conservative integrators preserve the mass of the solution up to machine precision.}\label{fig_planesource}
\end{figure}

\subsection{Planesource test case}\label{sec:planesource}
In this numerical example, the radiation transport equation \eqref{eq:RTE} in a one-dimensional slab geometry is solved. The computational domain is $x \in [-2,2]$ and $v \in [-1,1]$. The initial condition is a Gaussian with a standard deviation of $\sigma_{\mathrm{IC}} = 3\cdot 10^{-2}$ and is independent of $\mu$, meaning that particles are equally likely to have any particular direction. In slab geometry, the radiation transport equation \eqref{eq:RTE} then becomes
\begin{align*}
    \partial_t f + v \partial_x f + \sigma_s f =&\, \frac{\sigma_s}{2} \int_{-1}^1 f \, dv\, , \\
    f(t=0,x,v) =&\, \frac{1}{\sqrt{2\pi}\sigma_{\mathrm{IC}}}\exp \left( -\frac{x^2}{2\sigma^2} \right)\;.
\end{align*}
Here, we use isotropic scattering, that is the scattering kernel is independent of $v$ and we choose $k(x,v',v) = \sigma_s = 1$ as well as $\sigma_a = 0$. This test case is known as the planesource benchmark test \cite{ganapol2008analytical} and is widely used since an analytic solution to this problem is available. Furthermore, due to the moderate scattering strength and the fact that particles move in all directions, this test case is challenging for numerical methods and is often used to demonstrate deficiencies in a given numerical schemes. The numerical discretization is performed via an upwind flux using $1500$ spatial cells and $500$ moments with a CFL number of $0.9$, i.e., $\Delta t = 0.9\Delta x$.

We compare the conservative integrators BUG+Euler, Parallel+Euler, 2nd order Parallel+RK4, midpoint BUG+Euler, midpoint BUG+Heun, midpoint BUG+RK4 to the non-conservative BUG+Euler (i.e.~without adding $U(v)$ to the basis) integrator as baseline, see \Cref{tab_integrators} for an overview. 
We observe in \Cref{fig_planesource} that  all methods approximate the exact solution sufficiently well. Second order integrators oscillate less, but smear out slightly more. This is explained by the fact that the numerical flux is discretized with an upwind scheme that introduces significant numerical dissipation. The rank dynamics of the integrators are similar; the second order integrators exhibit slightly lower ranks.   
All conservative methods preserve the mass of the solution up to machine precision, whereas the non-conservative BUG integrator exhibits an  $10$ order of magnitude larger mass error than the conservative schemes. 


\subsection{Linesource test case}\label{sec:linesource}
Next, we examine the radiation transport equation in the context of the line source benchmark \cite{ganapol1999homogeneous}. This scenario involves a particle pulse distributed uniformly in all directions along the z-axis. The initial particle pulse is represented by a Gaussian distribution with a standard deviation of $\sigma_{\mathrm{IC}} = 0.03$. As time progresses, the particles interact with a background medium that has an isotropic scattering kernel $k(x,v,v') = \frac{\sigma_s}{4\pi}$ with $\sigma_s = 1$. Absorption is not present, that is, $\sigma_a = 0$. Then, for $v\in\mathbb{S}^2$ and $x\in\mathbb{R}^2$ the radiation transport equation reads
\begin{align*}
    \partial_t f + v_1 \partial_{x_1} f + v_2 \partial_{x_2} f + \sigma_s f =&\, \frac{\sigma_s}{4\pi} \int_{\mathbb{S}^2} f \, dv\, , \\
    f(t = 0,x,v) =&\, \frac{1}{4\pi\sigma_{\mathrm{IC}}^2} \exp\left(-\frac{\Vert x\Vert^2}{4\sigma_{\mathrm{IC}}^2}\right)\,.
\end{align*}

\begin{figure}[th]
    \centering
    \begin{subfigure}{0.32\textwidth}
        \includegraphics[width=\textwidth]{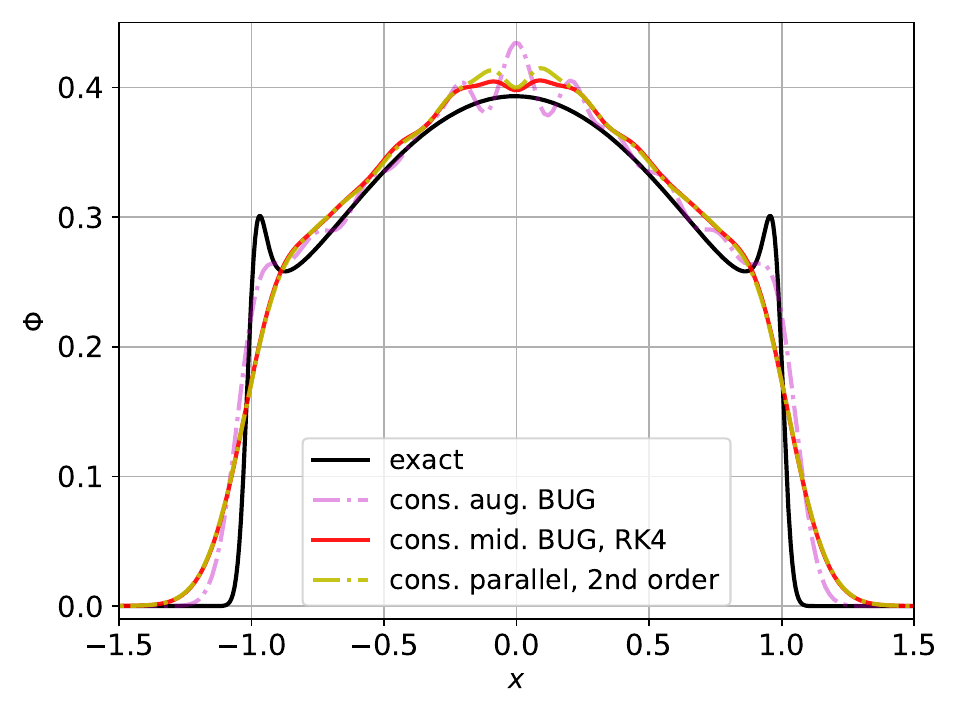}
        \caption{Horizontal cut of the solutions}
    \end{subfigure}
    \begin{subfigure}{0.32\textwidth}
        \includegraphics[width=\textwidth]{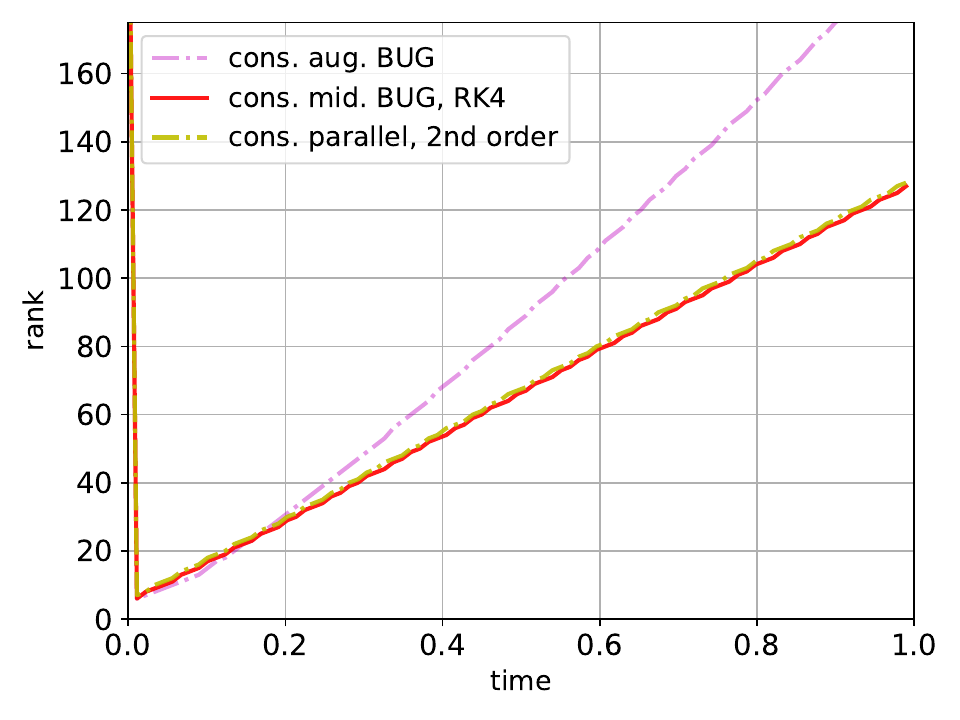}
        \caption{Rank evolution}
    \end{subfigure}
     \begin{subfigure}{0.32\textwidth}
        \includegraphics[width=\textwidth]{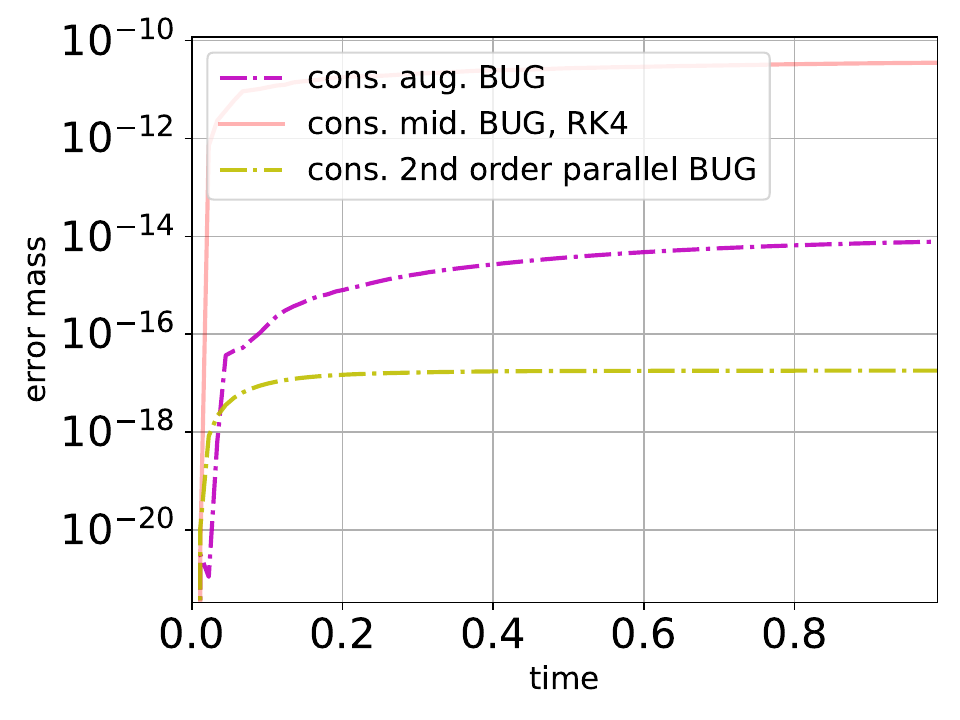}
        \caption{Mass error}
    \end{subfigure}
    \caption{2D2V Linesource. Solution mass across a horizontal cut through the domain (a)  rank evolution (b) and (c) mass error of the conservative integrators. Higher order integrators have a slower increasing rank, but smear out slightly more compared to first order integrators. \label{fig:scalar_flux_DLRA_adBUG_LineSource}}
\end{figure}

\begin{figure}[t!]
    \centering
    \begin{subfigure}{0.245\textwidth}
        \includegraphics[width=\textwidth]{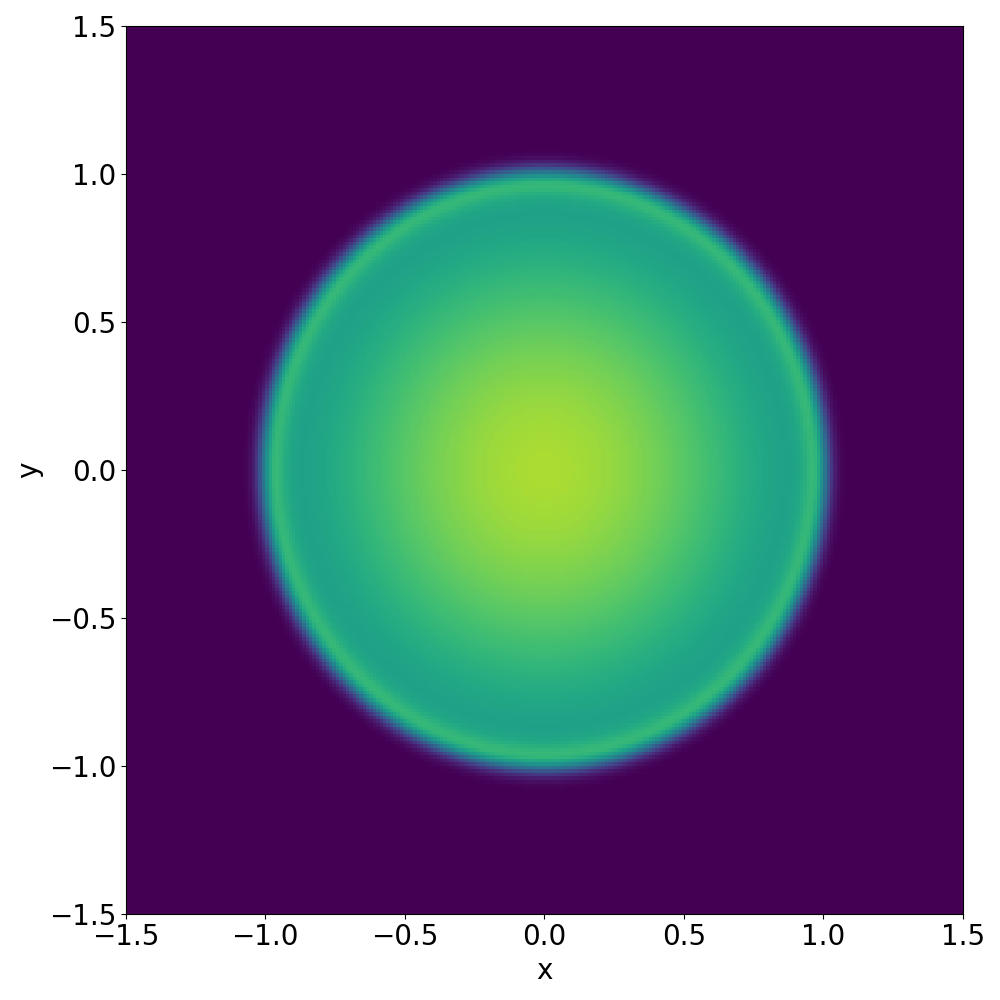}
        \caption{Analytic reference solution}
    \end{subfigure}%
    \begin{subfigure}{0.245\textwidth}
        \includegraphics[width=\textwidth]{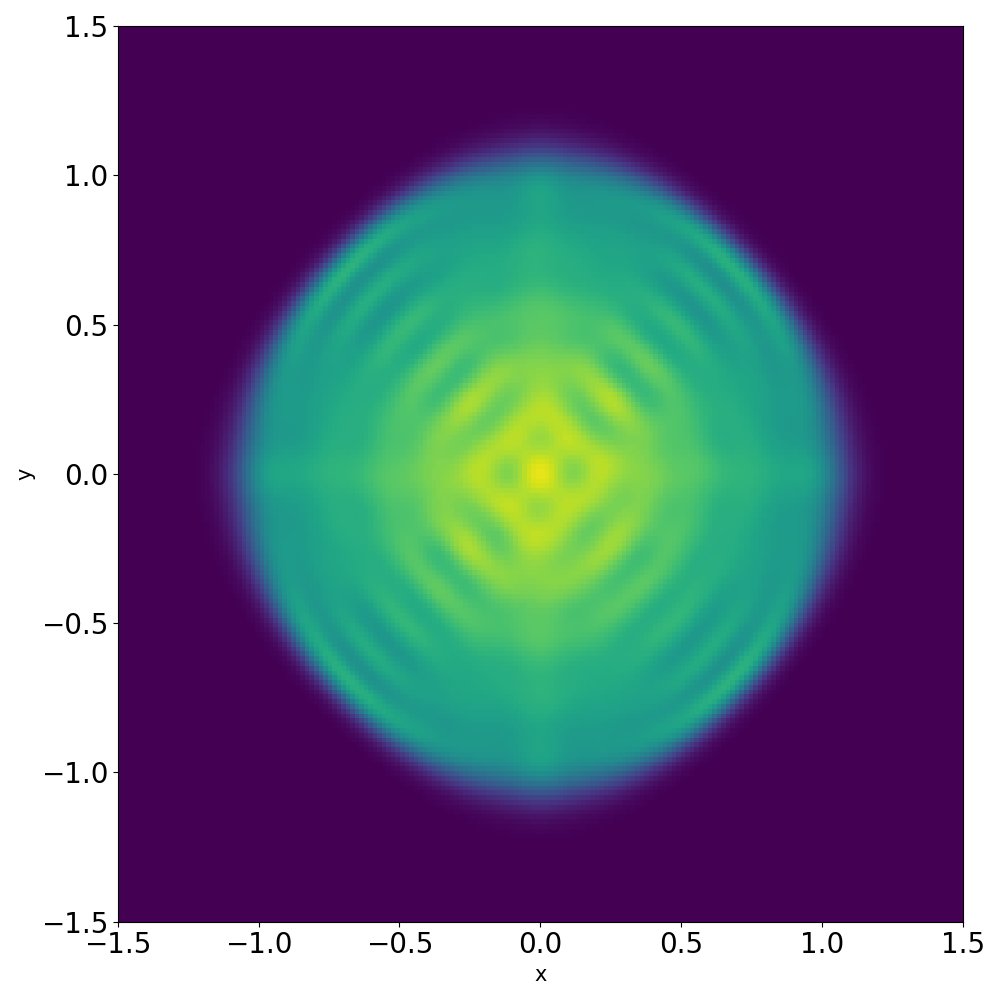}
        \caption{Cons. augmented BUG}
    \end{subfigure}%
    \begin{subfigure}{0.245\textwidth}
        \includegraphics[width=\textwidth]{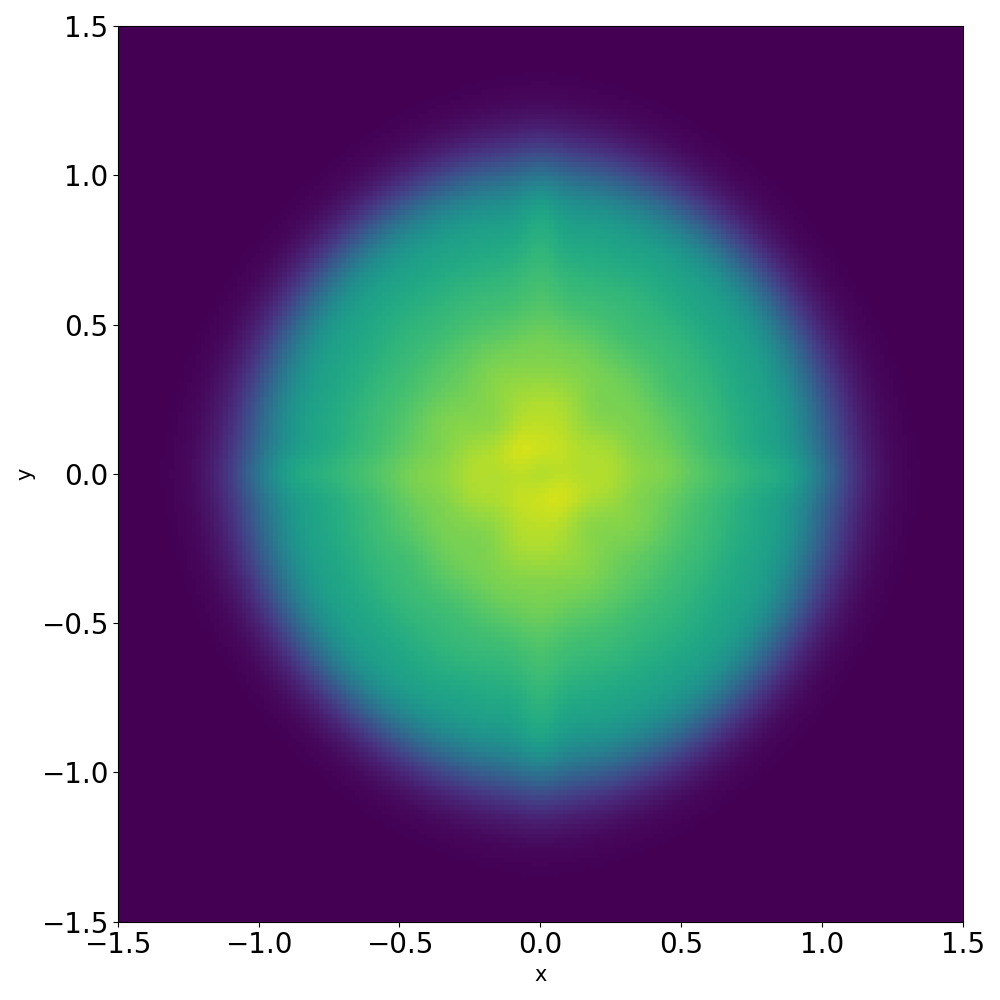}
        \caption{Cons. parallel  second order }
    \end{subfigure}
    \begin{subfigure}{0.245\textwidth}
        \includegraphics[width=\textwidth]{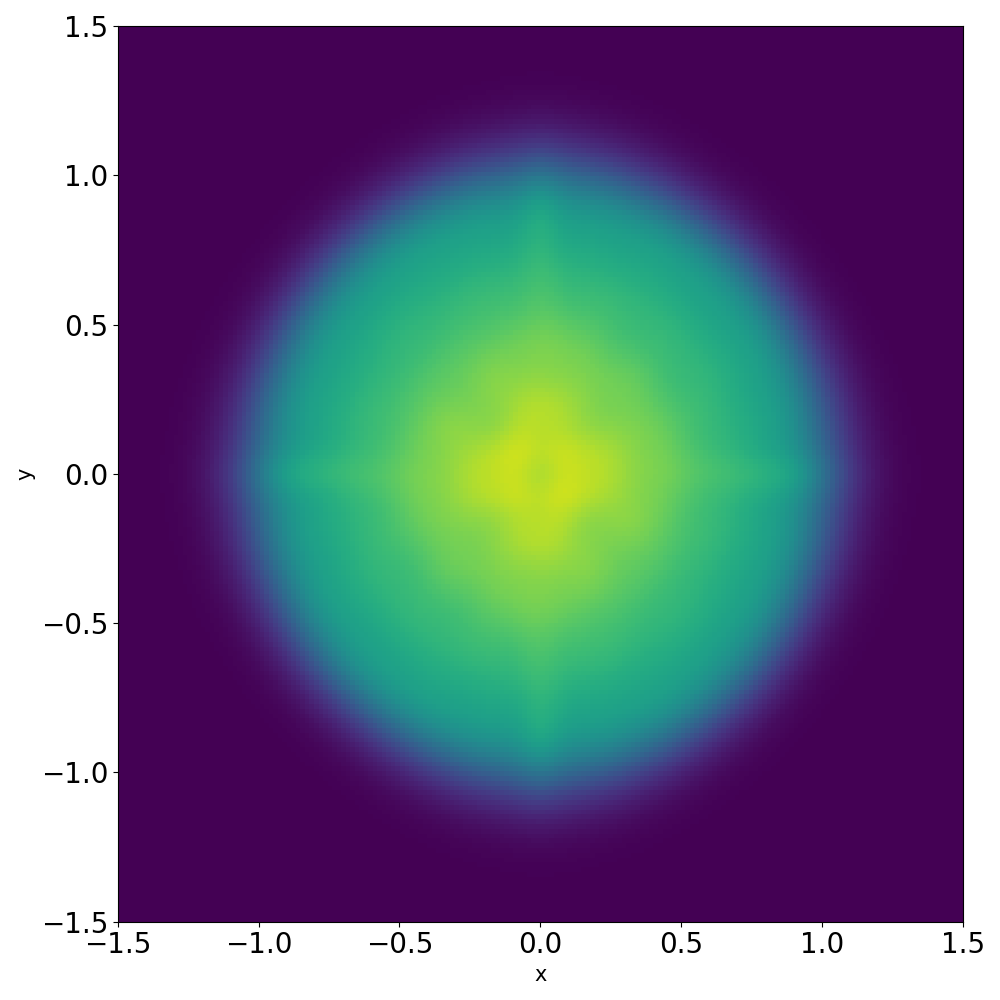}
        \caption{Cons. midpoint BUG, RK4 }
        \label{fig:scalar_flux_DLRA_fLukas_LineSource}
    \end{subfigure}
    \caption{\green{The scalar flux for the 2D2V Linesource simulation using the conservative augmented BUG integrator (b),  the conservative parallel second order BUG integrator with RK4 substeps (c),  and the conservative midpoint BUG integrator with RK4 substeps (d). The truncation threshold is given by $\vartheta=\bar{\vartheta} \Vert \widehat S^1 \Vert$ with $\bar\vartheta = 5\cdot 10^{-2}$, $t_{\mathrm{end}}= 1$.} \label{fig:scalar_flux_reference}}
\end{figure}

The scalar flux $\phi(t,x) := \int_{\mathbb{S}^2}f(t,x,v)\,dv$ of the analytic solution to this equation at time $t=1$ is depicted in Figure~\ref{fig:scalar_flux_reference}. Here, particles have moved into all possible directions on the sphere and are projected onto the two-dimensional spatial domain. 

In our numerical simulation, we use $250^2$ spatial cells to represent the spatial domain and $31^2$ modal coefficients to represent the directional basis. The chosen tolerance for the truncation step of the integrators is set to  $\vartheta=\bar{\vartheta} \Vert \widehat S^1 \Vert$ with $\bar\vartheta = 5\cdot 10^{-2}$. Similar to conventional numerical methods, dynamical low-rank approximation struggles to approximate this test case efficiently, since a large rank is required. The resulting scalar flux for the augmented BUG integrator with adaptive rank is depicted in Figure~\ref{fig:scalar_flux_DLRA_adBUG_LineSource}. We see oscillatory artifacts which are expected when using $31^2$ directional modes to represent the basis. 
It is apparent that higher order integrators produce less oscillations, but smear out the solution slightly more than first order low-rank integrators.

 In general, considering the challenging structure of this test case, the approximation results show a satisfactory solution quality. This solution quality is however reached with significantly increased ranks, which are depicted in Figure~\ref{fig:scalar_flux_DLRA_adBUG_LineSource}. 
 The {conservative augmented } BUG integrator for both the conservative system \eqref{eq:modifiedDLRA} and the original system \eqref{eq:ill-cond-dlra} yield conservation up to machine precision. 

\section{Conservation for the Vlasov--Poisson equations}\label{sec:discussion_vlasov}

For the radiative transfer equation, the schemes described in sections \ref{sec:local-cons-rbug} and \ref{sec:cons-parallel} can be applied directly. However, for the Vlasov--Poisson equations a number of additional steps, which we will describe in this section, need to be taken. The first complication arises because the velocity domain of the Vlasov--Poisson equation, in principle, is infinite. Thus, neither $V_1 \propto 1$ nor $V_2 \propto v$ lies in $L^2(\mathbb{R})$. Note that, for the sake of simplicity, we only consider a one-dimensional problem and mass and momentum conservation here. However, the extension to multiple dimensions is immediate, simply by augmenting with $v_1$, $v_2$, and $v_3$, separately. Of course, in any numerical scheme, the infinite domain is truncated. However, especially for $v$ this can still cause problems as it is not compatible with the (artificial) boundary conditions imposed in velocity space. For more details, we refer to \cite{einkemmer2021mass}. To overcome this we can, as has been suggested in \cite{einkemmer2021mass}, use a weighted space $L^2(\Omega_v, f_{0v})$. The low-rank approximation then takes the form (we follow the notation from \cite{einkemmer2021mass})
\begin{equation}
f(t,x,v)=f_{0v}(v) \sum_{i,j=1}^r X_i(t,x) S_{ij}(t) V_j(t,v), \label{eq:lr-weights}
\end{equation}
where $V_j \in L^2(\Omega, f_{0v})$ and $f_{0v}$ is a weight function. We will later discuss appropriate choices of $f_{0v}$. The main idea, however, is that it decays sufficiently rapidly such that $v$ lies in $L^2(\Omega, f_{0v})$.

For the Vlasov--Poisson equations, the dynamical low-rank equations of motions are as follows. \newline
\textit{K step:} 
\begin{equation}
\partial_t{K}_{j}(t,x)=-\sum_{l}c_{jl}^{1}\cdot\nabla_{x}K_{l}(t,x)+\sum_{l}c_{jl}^{2}\cdot E(K)(t,x)K_{l}(t,x)\label{eq:evolution-K}
\end{equation}
with
\[
c_{jl}^{1}=\langle V_j, V_l \rangle_v := \int_{\Omega_v} f_{0v} vV_{j} V_{l} \,\mathrm{d}v,\qquad 
c_{jl}^{2}= \int_{\Omega_v} V_{j} \nabla_{v}(f_{0v} V_{l})\,\mathrm{d}v.
\]
\textit{L step:} 
\begin{equation}
\partial_t{L}_{i}(t,v)  =\sum_{k}d_{ik}^{1}[E(L(t,\cdot))]\cdot\nabla_{v}(f_{0v}(v) L_{k}(t,v))-\sum_{k}(d_{ik}^{2}\cdot v)L_{k}(t,v),\label{eq:evolution-L}
\end{equation}
with
\[
d_{ik}^{1}[E]=\int_{\Omega_x} X_{i} E X_{k}\,\mathrm{d}x,\qquad\quad d_{ik}^{2}=\int_{\Omega_x} X_{i}(\nabla_{x}X_{k})\,\mathrm{d}x,
\]

\textit{S step:} 
\begin{equation}
\partial_t S_{ij}(t)  =\sum_{k,l}\left(c_{jl}^{1}\cdot d_{ik}^{2}-c_{jl}^{2}\cdot d_{ik}^{1}[E(S(t))]\right)S_{kl}(t)\label{eq:evolution-S}
\end{equation}

Note that they differ from \cite{einkemmer2018low} by the fact that the inner products and some of the terms include the weight $f_{0v}$, as we use the modified low-rank expression given in equation \eqref{eq:lr-weights}. We also note that we have not done the modification in \cite{einkemmer2021mass,einkemmer2023robust} to make the scheme conservative on the continuous level as, in contrast to the methods discussed in that paper, this procedure is not necessary to obtain a conservative scheme on the discrete level for the schemes discussed here. It is interesting to note in that context that, as mentioned in section \ref{sec:local-cons-rbug}, it is in general not possible for the approach described here to be conservative on the continuous level. However, once the time discretization has been performed, we obtain a conservative scheme. This is also the reason why the augmentation is specific to the time integration method chosen.

In addition to the introduction of the weight, a number of additional considerations have to be taken into account in order to obtain a robust and efficient scheme that is mass and momentum conservative up to machine precision. We will detail them in the following subsections.

\subsection{Approximation of the electric field \label{sec:approx-E}}

If we naively apply our conservative scheme to the equations of motion \eqref{eq:evolution-K}-\eqref{eq:evolution-S}, the electric field is recomputed in each of the three steps and in each stage of the Runge--Kutta method used for the corresponding time integration. Since evaluating the electric field requires a Poisson solve, this is expensive and is rarely done in practice. Usually, first order methods simply compute the electric field once at the beginning of the time step and use that value for the entire computation. For second order methods, it is sufficient to obtain a first-order approximation at time $t^{n+1/2}$. Thus, two Poisson solver are sufficient. The first to compute $E^{n}$ which is used to obtain an approximation of $f$ at $t^{n+1/2}$ and the second to compute $E^{n+1/2}$ from that value. Dynamical low-rank methods for the Vlasov--Poisson equation have used similar techniques to keep the computational cost low (see, e.g., the algorithms presented in \cite{einkemmer2018low,cassini2022efficient,guo2022conservative}). If, on the other hand, this is done naively, we require $3s$ Poisson solves, where $s$ is the number of stages of the Runge--Kutta method.

However, if we desire a scheme that conserves momentum exactly using $E^{n+1/2}$ in the algorithm is not sufficient. This can be easily seen by considering (for simplicity in 1+1d)
\[ \partial_t f + v \partial_x f(t,x,v) - E^{n+1/2}  \partial_v f(t,x,v) = 0. \]
Multiplying by $v$ and integrating in velocity yields
\[ \partial_t j + \partial_x  \int v f(t,x,v) \,\mathrm{d}v = E^{n+1/2} \int \partial_v f(t,x,v) v \,\mathrm{d}v. \]
The term on the right-hand side gives
\[ E^{n+1/2} \int v \partial_v f(t,x,v) \,\mathrm{d}v = -E^{n+1/2} \rho(t,x), \]
but this is not in conservative form as
\[ \int E^{n+1/2} \rho(t,x) \,\mathrm{d}x = - \int E^{n+1/2} \partial_x E(t,x)  \,\mathrm{d}x.  \]
If instead of $E^{n+1/2}$ the correct field $E(t,x)$ is used, the term can easily be written in conservative form as $\tfrac{1}{2}\partial_x (E^2)$, but this is not possible here.

Thus, such an approximation, even on the continuous level, violates momentum conservation. The question then is, can we obtain a momentum conservative scheme while still keeping the number of Poisson solves to a minimum? In this context, it is important to realize that for showing conservation only the final S step of the scheme is important. The remainder only acts to build a sufficiently accurate basis. Thus, for example, for the BUG and midpoint BUG integrators computing everything up to the final S step with $E^n$ and $E^{n+1/2}$ is sufficient to obtain a first and second order scheme, respectively. In the final S step one Poisson solve for each stage of the used Runge--Kutta method is required. This reduces the total number of Poisson solves from $3s$ to $2+s$, which is a significant improvement.

The conservative parallel integrator of order two is an exception. Since in this case the K step is used to reconstruct $S$, we have to compute the appropriate electric field in each stage of the Runge--Kutta method used for the K step to obtain conservation. In the L step we can still use the cheaper approximation of $E^{n+1/2}$. 

\subsection{Appropriate choice of the weight function}

In principle, the weight function $f_{0v}$ can be chosen arbitrarily as long as it decays sufficiently fast to ensure that $v \in L^2(\Omega_v, f_{0v})$. For example, in \cite{einkemmer2021mass} it is suggested to use a Maxwellian, i.e. $\exp(-v^2/2)$. While technically this choice is fine, it can result in inaccurate results, especially in the nonlinear regime. The problem is that the projection in the dynamical low-rank algorithm minimizes the error in the given norm. This, in particular, implies that the error for large velocities, where $\exp(-v^2/2)$, is small, can become relatively large. 

\begin{figure}[t]
\begin{center}\includegraphics[width=0.5\textwidth]{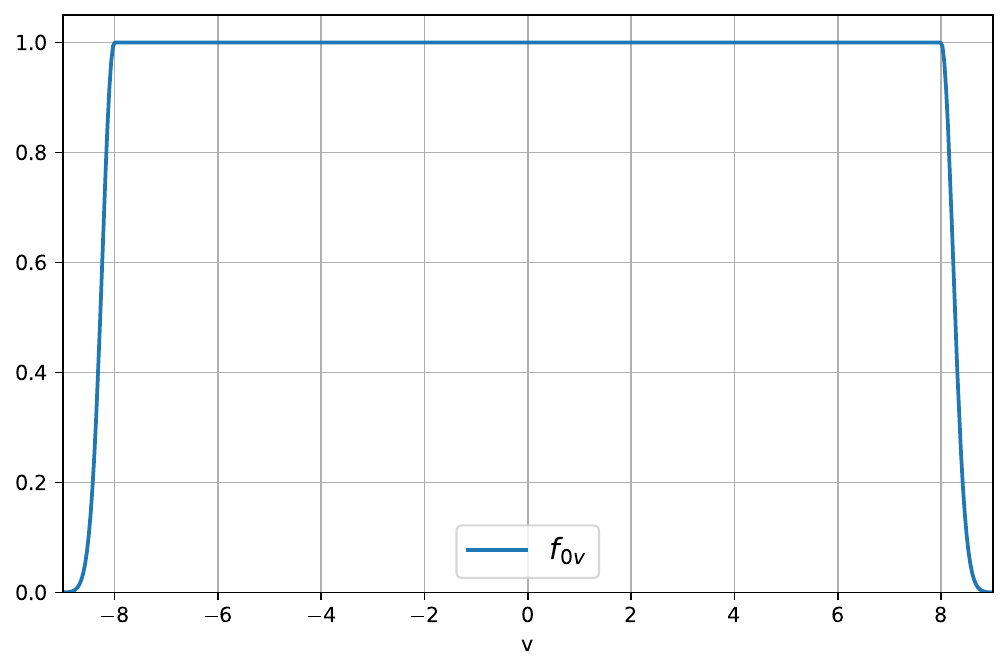}\end{center}
\caption{Plot of equation \eqref{eq:f0v} with $M=10^{-4}$, $v_0=8$, and $v_1=9$. \label{fig:f0v}}
\end{figure}

On the other hand, the weight function on the boundary should be small such that $v f_{0v}$ still satisfies the imposed artificial boundary conditions. We have found that in practice
\begin{equation} f_{0v}(v) = \begin{cases}
1 & v\leq v_{0}\land v\geq-v_{0}\\
\exp\left(\log(M)\left(\frac{v-v_{0}}{v_{1}-v_{0}}\right)^{2}\right) & v>v_{0}\\
\exp\left(\log(M)\left(\frac{v+v_{0}}{v_{1}-v_{0}}\right)^{2}\right) & v<v_{0}
\end{cases}, \label{eq:f0v} \end{equation}
with $M=10^{-4}$, $v_1$ the velocity at the boundary of the velocity domain, and $v_0$ a velocity close to the boundary, gives good results. A plot of the this weight function is shown in Figure \ref{fig:f0v} and we have used those parameters in all our numerical simulations of the Vlasov--Poisson equations in section \ref{sec:vp-numerics}.

\subsection{Avoiding unfavorable round-off errors}

We have shown that the schemes described are both locally and globally conservative. However, in practical computation, there is one additional source of error; namely, the finite precision representation of real numbers. In the analysis we have neglected this error, as is done almost universally in the literature. The justification for this is that this round-off error is much smaller (approximately $10^{-16}$ for double precision floating point numbers) compared to the numerical error introduced by the dynamical low-rank approximation and the discretization of the problem. However, if one implements the algorithm in a straightforward way, it can still happen that round-off errors propagate in an unfavorable way which, after some time, gives an error in the conserved quantities that is much larger than machine precision. In this section, we will explain why this is a problem for conservative dynamical low-rank methods and how to implement the numerical method such as to avoid unfavorable propagation of round-off errors.

In order to implement the K and S step we have to compute
\[ c^2_{jl} = \int V_j \nabla_v (f_{0v} V_l) \,dx. \]
In particular, we note that due to our augmentation $V_1 \propto 1$ and $V_2 \propto v$, we get by integration by parts and orthogonality
\begin{equation} c^2_{1,l} = 0 \text{ for all } l 
\qquad \text{and} \qquad
\qquad c^2_{2,l} = 0 \text{ for all } l\geq 2. \label{eq:c2-rel} 
\end{equation}

After performing a space discretization $c^2$ is approximated by a quadrature rule, e.g. the trapezoidal rule. Then
\begin{equation} \label{eq:discrete-c2}
    c^2 \approx h_v V^{\top} D_v \text{diag}(\{f_{0v}(v_k)\}_k) V,
\end{equation}
where $V \in \mathbb{R}^{n_v \times r}$ are the discretized low-rank factors, $D_v$ is the discretized differential operator, and $v_k$ are the grid points. Assuming that $D_x$ is chosen such that integration by parts (or rather summation by parts) holds true, which e.g.~is the case for standard centered differences, we obtain the same relations as given in equation \eqref{eq:c2-rel}. 

However, numerically we observe that equations \eqref{eq:c2-rel} do not hold true and that the error is commonly orders of magnitude larger than machine precision. The reason for this is that for problems where some singular values of $S$ are small (which is commonly the case), not all $V$ are well defined. In this case columns of $L$ are close to machine precision and the QR decomposition has thus insufficient data to obtain corresponding rows in $V$ from the problem. What then commonly happens is that whatever algorithm is used for the QR decomposition generates orthogonal columns that are essentially random. For the algorithm, this does not matter since once $V$ is multiplied by $S$ the result is  zero no matter what the value in these columns of $V$ are. However, these columns are commonly very rough, i.e.~are not discretizations of smooth functions. Then applying the differential operator $D_x$ generates columns with huge entries and this raises the error in $c^2$ to orders of magnitude above machine precision. This has a detrimental effect on conservation as our scheme is only conservative up to the error in $c^2$. In the implementation, we thus have to explicitly enforce equation \eqref{eq:c2-rel}. We do this simply by setting the corresponding entries to zero after an approximation has been obtained by evaluating equation \eqref{eq:discrete-c2}.

Another point is that it is generaly more favorable to compute a discretization of
\[ \partial_v (f_{0v}(v)V_l(v)) \]
as opposed to a discretization
\[ V_l \partial_v f_{0v} + f_{0v} \partial_v V_l(v). \]
Although $\partial_v f_{0v}$ is often available analytically, $\partial_v V_l(v)$ does not necessarily satisfy the correct boundary conditions. For example, for $V_2 \propto v$ imposing periodic or Dirichlet boundary conditions creates large errors on the boundary when applying $D_v$. On the other hand, this is not problematic for $f_{0v} V$ since $f_{0v}$ is chosen to be small at the boundary.  

We also emphasize here that orthogonality as close as possible to machine precision is crucial to obtain good conservation over long times. In particular, in our view, it is good practice to explicitly perform an orthonormalization of $V_1$ and $V_2$ even if they have been chosen to be exactly orthogonal. Both numerical and round-off errors can contribute to the different columns of $V$ to be not exactly orthogonal in the actual computation, which can degrade conservation. 

\begin{figure}[t]
\includegraphics[width=\textwidth]{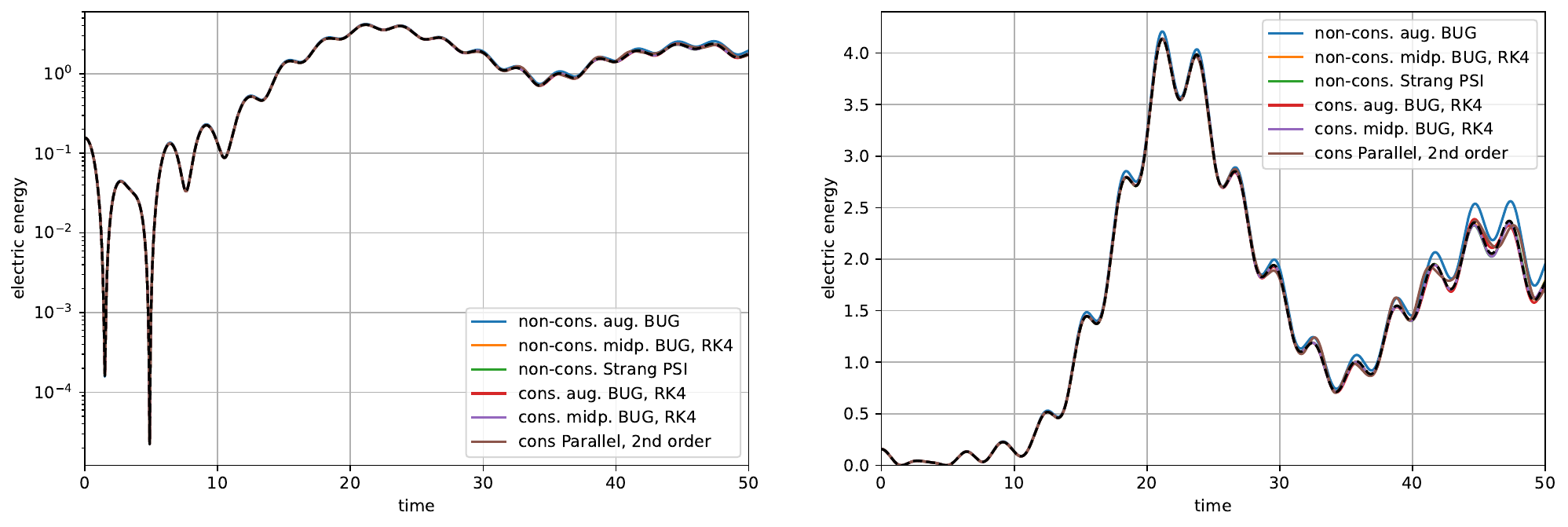}
\includegraphics[width=\textwidth]{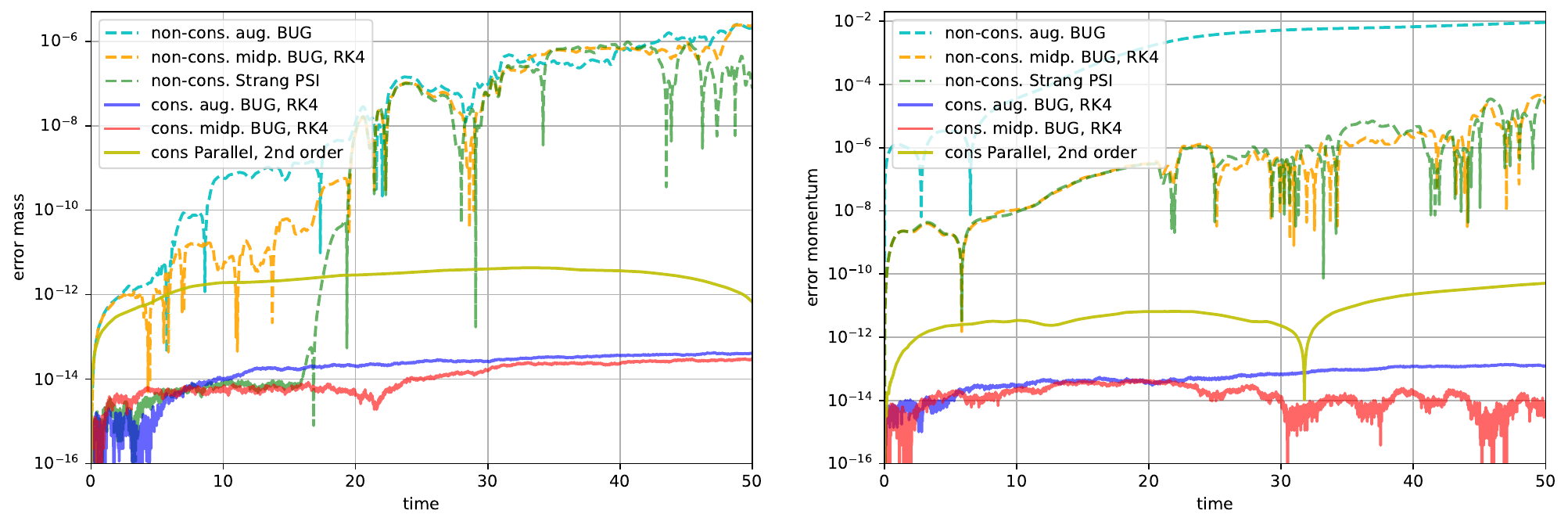}
\caption{Time evolution of the electric energy for a bump-on-tail instability is shown on the top. On the bottom the time evolution of the relative error in mass (left) and momentum (right) is shown. \label{fig:vp-evol-cons}}
\end{figure}

\section{Numerical results for the Vlasov--Poisson equations \label{sec:vp-numerics}}

As an example, we solve the 1+1d Vlasov--Poisson equations for a bump-on-tail initial condition given by
\[ f(0,x,v) = (1.0+ 0.03 \cos(0.3 x)) \left(0.9 \frac{\exp(-v^2/2)}{\sqrt{2 \pi}} + 0.1 \frac{\exp(-2(v-4.5)^2)}{\sqrt{\pi/2}} \right), \]
where $\Omega_x = [0, 20\pi]$ and $\Omega_v = [-9, 9]$. In all simulations, we use the standard centered difference scheme of order two (which in particular satisfies the summation by parts property) and 128 and grid points in both directions. As time integrator we choose the classic method of Runge and Kutta (RK4), which is stable in combination with a central difference scheme as the domain of stability includes part of the imaginary axis. For the simulations we use $r=25$, for non-conservative methods, and $r=27$ for conservative methods (in order to have the same number of independent basis functions in both cases) and a time step size of $\Delta t = 0.01$.  For the reference simulation we use a Strang projector splitting integrator with $r=100$ and time step size $\Delta t = 10^{-4}$.

The time evolution of the electric energy and the error committed in mass and momentum for a number of dynamical low-rank schemes are shown in Figure \ref{fig:vp-evol-cons}. As expected, the classic BUG, midpoint BUG and Strang projector splitting commit an error that is far above machine precision. Especially for momentum conservation the second order methods (midpoint BUG and Strang) perform better but still have an error above $10^{-5}$. The conservative integrators proposed in this work (here exemplified by mod BUG and mod BUG midp) show conservation up to machine precision for both mass and momentum.

\begin{figure}[t]
\begin{center}\includegraphics[width=0.5\textwidth]{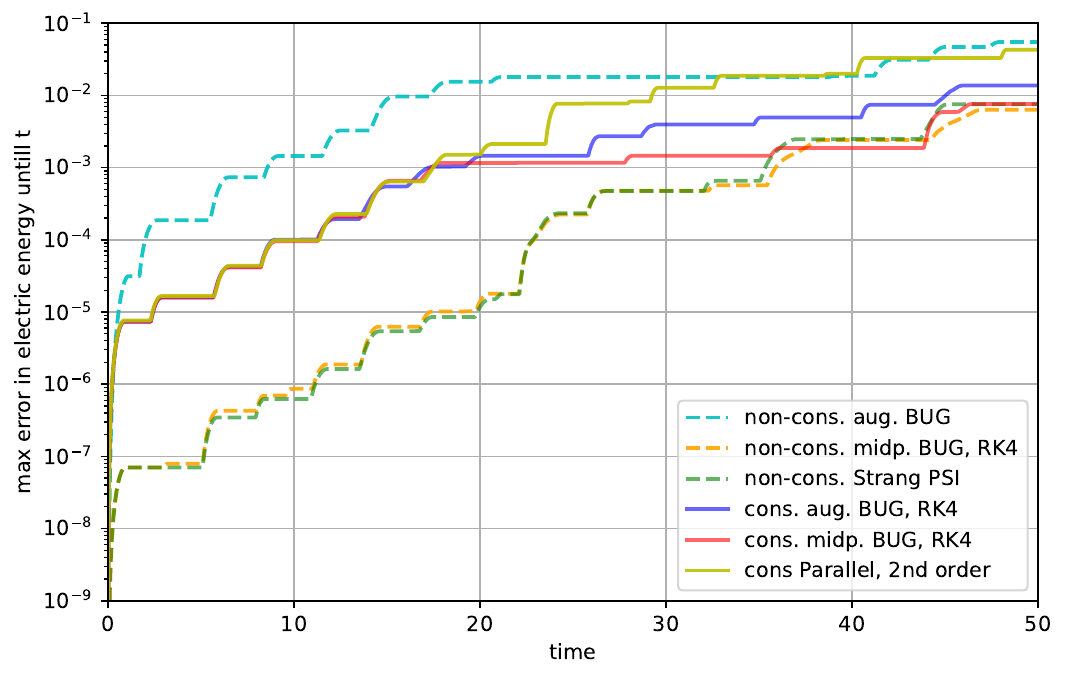}\end{center}
\caption{Time evolution of the error in the electric field for a bump-on-tail instability. The maximal error up to time $t$ is shown in this plot.\label{fig:vp-error}}
\end{figure}

We note that except perhaps for the BUG integrator all methods match the reference solution very well. In Figure \ref{fig:vp-error} we investigate this in more detail by plotting the error in the electric energy up to time $t$. The latter is done to remove oscillations in the error and provide a more legible plot. We see that the error in the classic BUG integrator is relatively large, presumably since it only uses a first order approximation of the electric field. The modified BUG integrator reduces the error significantly. Note, however, that in this case we compute the electric field in each stage of RK4 method in the S step, see section \ref{sec:approx-E} for more details. Thus, the integrator is more accurate but also more expensive. The second order methods are more accurate still and the overall accuracy with respect to the electric field is comparable for conservative and non-conservative methods, as we would expect.

\begin{figure}[t]
\includegraphics[width=\textwidth]{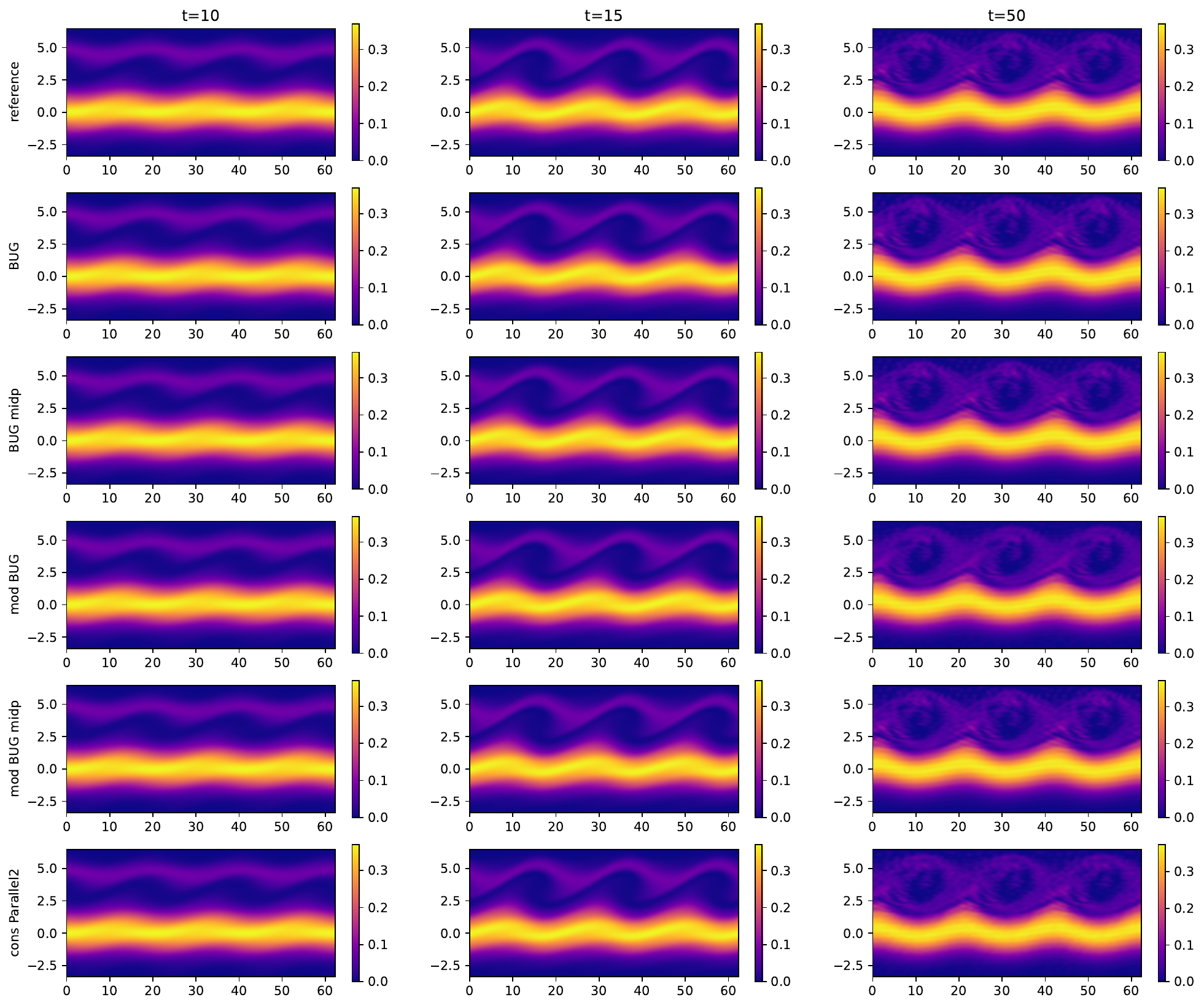}
\caption{Distribution function for the bump-on-tail instability. We plot different methods along the rows and different times along the columns.\label{fig:vp-distribution}}
\end{figure}

We also have plotted the distribution function $f$ for various integrators at different times (see Figure \ref{fig:vp-distribution}). We observe that they all match well compared to the reference solution and also with what has been reported in the literature.

\section{Conclusion}
In the work, we have derived conservative extensions to BUG integrators and demonstrated their effectiveness for two kinetic problems. The construction has been based on the observation that a conservative scheme must preserve conservative basis functions (a property that is achieved through a conservative truncation) and the exact representation of local conservation terms by the spatial basis (which can be ensured by an additional augmentation step). We have demonstrated three cases where such an augmentation step is not required, one being the parallel integrator, which allows for conservation with minor modifications and independent of the specific time integrator. We provide details on the implementation for both radiative transfer and plasma physics applications and have pointed out sources of errors that can destroy conservation.

\section{Acknowledgements}
The work of Steffen Schotthöfer is sponsored by the Office of Advanced Scientific Computing Research, U.S. Department of Energy, and performed at the Oak Ridge National Laboratory, which is managed by UT-Battelle, LLC under Contract No. DE-AC05-00OR22725 with the U.S. Department of Energy. The United States Government retains and the publisher, by accepting the article for publication, acknowledges that the United States Government retains a non-exclusive, paid-up, irrevocable, world-wide license to publish or reproduce the published form of this manuscript, or allow others to do so, for United States Government purposes. The Department of Energy will provide public access to these results of federally sponsored research in accordance with the DOE Public Access Plan (http://energy.gov/downloads/doe-public-access-plan).
\bibliographystyle{abbrv}
\bibliography{main} 
	
\end{document}